\documentclass[a4paper,12pt]{article}

\newtheorem{theo}{Theorem}[section]
\newtheorem{ex}[theo]{Example}
\newtheorem{prop}[theo]{Proposition}
\newtheorem{lem}[theo]{Lemma}
\newtheorem{cor}[theo]{Corollary}

\newtheorem{rema}[theo]{Remark}
\newtheorem{remas}[theo]{Remarks}
\newtheorem{propbis}{Proposition}

\def \kbar {{\bar k}}

\def \dem {\paragraph{Proof : }}

\def \Romannumeral #1 {\expandafter\uppercase\expandafter {\romannumeral #1} }
\def \br {{\rm{Br\,}}}

\def \Si {{\Sigma}}

\def \P {{\bf P}}
\def \tors {{\rm tors}}
\def \Ga {{\Gamma}}

\def \pic {{\rm {Pic\,}}}

\def \gal {{\rm{Gal\,}}}

\def \calo {{\cal O}}
\def \spec {{\rm{Spec\,}}}

\def \Hom {{\rm {Hom}}}

\def\ra{\rightarrow}

\def \nr {{\rm nr}}
\def\k{\kappa}
\def \Z {{\bf Z}}
\def \Q {{\bf Q}}
\def \F {{\bf F}}

\def \RR {{\bf R}}
\def \C {{\bf C}}

\def \cok {{\rm{coker\,}}}
\def \im {{\rm {Im\,}}}
\def \G {{\bf G}_m}

\def \mathopto {\mathop{\to}}
\def\smallsquare{\vbox{\hrule\hbox{\vrule height 1 ex\kern 1 ex\vrule}\hrule}}
\def\enddem{\hfill \smallsquare\vskip 3mm}

\def \abstract{\paragraph{Abstract. }}
\usepackage{amscd}
\usepackage{amssymb}
\usepackage{amsmath}
\everymath{\displaystyle}

\title{Arithmetic duality theorems for 1-motives}
\author{David Harari and Tam\'as Szamuely}
\date{}

\def \nr {{\rm nr}}
\def \ext {{\rm Ext}^1}

\def \tors {{\rm tors}}
\def \id {{\rm id}}

\def \fp {{\bf {\cal F}}}

\def \hyp {{\bf H}}
\def \thyp{\widehat{{\bf H}}}

\DeclareFontFamily{U}{wncy}{}
\DeclareFontShape{U}{wncy}{m}{n}{%
   <5>wncyr5%
   <6>wncyr6%
   <7>wncyr7%
   <8>wncyr8%
   <9>wncyr9%
   <10>wncyr10%
   <11>wncyr10%
   <12>wncyr6%
   <14>wncyr7%
   <17>wncyr8%
   <20>wncyr10%
   <25>wncyr10}{}
\DeclareMathAlphabet{\cyrille}{U}{wncy}{m}{n}
\def\Sha{\cyrille X}
\def \R{{\bf R}}

\begin{document}
\maketitle

\setcounter{section}{-1}

\section{Introduction}

Duality theorems for the Galois cohomology of commutative group schemes over local and global fields are among the most fundamental results in arithmetic. Let us briefly and informally recall some of the most famous ones.

Perhaps the earliest such result is the following. Given an algebraic torus $T$ with character group $Y^{\ast}$ defined over a $p$-adic field $K$, cup-products together with the isomorphism $\br(K)=H^2(K, \G)\cong \Q/\Z$ given by the invariant map of the Brauer group of $K$ define canonical pairings
$$
H^i(K, T)\times H^{2-i}(K, Y^{\ast})\to\Q/\Z
$$
for $i=0,1, 2$. The Tate-Nakayama duality theorem (whose original form can be found in \cite{tatenak}) then asserts that these pairings become perfect if in the cases $i\neq 1$ we replace the groups $H^0$ by their profinite completions. Note that this theorem subsumes the reciprocity isomorphism of local class field theory which is the case $i=0$, $T=\G$.

Next, in his influential expos\'e \cite{wc}, Tate observed that given an abelian variety $A$, the Poincar\'e pairing between $A$ and its dual $A^{\ast}$ enables one to construct similar pairings  
$$
H^i(K, A)\times H^{1-i}(K, A^{\ast})\to\Q/\Z
$$
for $i=0,1$ and he proved that these pairings are also perfect.

The last result we recall is also due to Tate. Consider now an abelian variety $A$ over a number field $k$, and denote by $\Sha^1(A)$ the Tate-Shafarevich group formed by isomorphism classes of torsors under $A$ that split over each completion of $k$. Then Tate constructed a duality pairing 
$$
\Sha^1(A)\times \Sha^{1}(A^{\ast})\to\Q/\Z
$$  
(generalising earlier work of Cassels on elliptic curves) and announced in \cite{stockholm} that this pairing is nondegenerate modulo divisible subgroups or else, if one believes the widely known conjecture on the finiteness of $\Sha^1(A)$, it is a perfect pairing of finite abelian groups. Similar results for tori are attributed to Kottwitz in the literature; indeed, the references \cite{ko1} and \cite{ko2} contain such statements, but without (complete) proofs.

In this paper we establish common generalisations of the results mentioned above for {\em 1-motives.} Recall that according to Deligne, a 1-motive over a field $F$ is a two-term complex $M$ of $F$-group schemes $[Y\to G]$ (placed in degrees -1 and 0), where $Y$ is the $F$-group scheme associated to a finitely generated free abelian group equipped with a continuous $\gal(F)$-action and $G$ is a semi-abelian variety over $F$, i.e. an extension of an abelian variety $A$ by a torus $T$. As we shall recall in the next section, every 1-motive $M$ 
as above has a {\em Cartier dual} $M^{\ast}=[Y^{\ast}\to G^{\ast}]$ equipped with a canonical (derived) pairing $M\otimes^{\bf L}M^{\ast}\to\G[1]$ generalising the ones used above in the cases $M=[0\to T]$ and $M=[0\to A]$. This enables one to construct duality pairings for the Galois hypercohomology groups of $M$ and $M^{\ast}$ over local and global fields.

Let us now state our main results. In Section 2, we shall prove:

\begin{theo}
Let $K$ be a local field and let $M=[Y \ra G]$ be a 1-motive over $K$. 
For $i=-1, 0, 1, 2$ there are canonical pairings
$$
\hyp^i(K, M)\times \hyp^{1-i}(K, M^{\ast})\to\Q/\Z
$$
inducing perfect pairings between
\begin{enumerate}
\item the profinite group $\hyp^{-1}_{\wedge}(K,M)$ and the discrete group 
$\hyp^2(K,M^{\ast})$;
\item the profinite group $\hyp^0(K,M)^{\wedge}$ and the discrete group 
$\hyp^1(K,M^{\ast})$.
\end{enumerate}
\end{theo}

Here the groups $\hyp^0(K,M)^{\wedge}$ and $\hyp^{-1}_{\wedge}(K,M)$ are obtained from the corresponding hypercohomology groups by certain completion procedures explained in Section~2. We shall also prove there a generalisation of the above theorem to 1-motives over henselian local fields of mixed characteristic and show that in the duality pairing the unramified parts of the cohomology are exact annihilators of each other.

Now let $M$ be a 1-motive over a number field $k$. For all $i\geq 0$ define the Tate-Shafarevich groups
$$
\Sha^i(M)={\rm Ker} \, [\hyp^i(k, M)\to\prod_v\hyp^i(\hat k_v, M)]
$$
where the product is taken over completions of $k$ at all (finite and infinite) places of $k$. Our main result can then be summarised as follows.

\begin{theo}\label{second}
Let $k$ be a number field and $M$ a 1-motive over $k$. There exist canonical 
pairings
$$
\Sha^i(M) \times \Sha^{2-i}(M^{\ast})\to\Q/\Z
$$
for $i=0, 1$ which are non-degenerate modulo divisible subgroups.
\end{theo}

In fact, here for $i=0$ the group $\Sha^0(M)$ is finite, 
so the left kernel is trivial. 
Assuming the finiteness of the (usual) Tate-Shafarevich group of an 
abelian variety one derives that for $i=1$ the pairing is a perfect pairing of finite groups.

The pairings used here can be defined purely in terms of Galois cohomology (see Section 6); however, to prove the duality isomorphisms
we first construct pairings using \'etale cohomology in Sections 3 and 4, and then in the last section we compare them to the Galois-cohomological one which in the case of abelian varieties gives back the classical construction of Tate.

Finally, in Section~5 we establish a twelve-term Poitou-Tate type exact sequence similar to the one for finite modules, assuming the finiteness of the Tate-Shafarevich group. The reader is invited to look up the precise statement there.

Since this is a paper containing the word ``motive", it is appropriate to explain our motivations for establishing the generalisations offered here. The first of these should be clear from the above: working in the context of 1-motives gives a unified and symmetric point of view on the classical duality theorems cited above and gives more complete results than those known before. As an example, we may cite the duality between $\Sha^i(T)$ and $\Sha^{3-i}(Y^{\ast})$ 
($i=1,2$) for an algebraic torus $T$ with character group $Y^{\ast}$ which is a special case of Theorem \ref{second} above (see Section~4); 
it is puzzling to note that the reference \cite{neukirch} only contains 
the case $i=1$, whereas \cite{adt} only the case $i=2$.
Another obvious reason is that due to recent spectacular progress in the theory of mixed motives there has been a regain of interest in 1-motives as well; indeed, the category of 1-motives over a field (with obvious morphisms) is equivalent, up to torsion, to the subcategory of the triangulated category of mixed motives (as defined, e.g., by Voevodsky) generated by motives of varieties of dimension at most 1. 

But there is a motivation coming solely from the arithmetic duality theory. In fact, if one tries to generalise the classical duality theorems of Tate to a semi-abelian variety $G$, one is already confronted to the fact that the only reasonable definition for the dual of $G$ is the dual $[Y^{\ast}\to A^{\ast}]$ of the 1-motive $[0\to G]$, where actually $A^{\ast}$ is the dual of the abelian quotient of $G$ and $Y^{\ast}$ is the character group of its toric part. Duality results of this type are needed for the study of the arithmetic of $G$ over local and global fields: for instance, such a question (which we shall address elsewhere) is the unicity of the Brauer-Manin obstruction to the Hasse principle for rational points on torsors under $G$. (This question was raised by Skorobogatov 
in \cite{skobook}, p. 133). 

In conclusion to this introduction, we  cannot resist the temptation of recalling the basic example of a dual pair of 1-motives ``appearing in nature''. For this, let $X$ be a smooth projective variety over an algebraically closed field. It is well-known that the Albanese variety ${\rm Alb}_X$ and the Picard 
variety $\pic^0_X$ of $X$ are dual abelian varieties. Now consider an open subvariety $U\subset X$. For such a $U$ a natural object to consider is the generalised Albanese variety $\widetilde{{\rm Alb}}_U$ defined by Serre \cite{moruniv}: it is a semi-abelian variety equipped with a morphism $U\to\widetilde{{\rm Alb}}_U$ universal for morphisms of $U$ into semi-abelian varieties that send a fixed base point to 0. According to \cite{moruniv2} (see also \cite{rama}) the dual of the 1-motive $[0\to \widetilde{{\rm Alb}}_U]$ is the 1-motive $[Y\to \pic^0_X]$, where $Y$ is the group of divisors on $X$ supported in $X\setminus U$ and algebraically equivalent to 0, and the map is the divisor class map.

After the first version of this paper has been posted on the Web,
N. Ramachandran has kindly drawn our attention to the paper of Berkovich
\cite{berkovich} where duality theorems for the Galois cohomology of
commutative group schemes over local and global fields are proven using the
technique of satellites. Though some of his results are certainly not
unrelated to ours, they seem to give the same statements only in the classical
case of abelian varieties. 

We thank Yves Andr\'e, Jean-Louis Colliot-Th\'el\`ene, Ofer Gabber, Joost van Hamel, Bruno Kahn, Niranjan Ramachandran and Michael Spiess for pertinent comments. The referee's remarks were very helpful for us. Needless to say, we owe a considerable intellectual
debt to the work of Tate, Artin/Verdier, Milne and others rendered in
\cite{adt}. A major part of this work was done during the second author's stay at the Universit\'e de Paris-Sud whose generous hospitality is gratefully acknowledged. The second author was also partially supported by OTKA grants Nos. T-031984 and F-032325.

\bigskip
\bigskip

{\bf \large Some notation and conventions}

\smallskip

Let $B$ be an abelian group.
For each integer $n >0$, $B[n]$ stands for the $n$-torsion subgroup of $B$
and $B_{\rm tors}$ for the whole torsion subgroup of $B$. We shall often
abbreviate the quotient $B/nB$ by $B/n$.
For any prime number $\ell$, we denote by $B\{\ell\}$ the
$\ell$-primary torsion subgroup of $B$ and by $\bar B\{\ell\}$ the quotient 
of $B\{\ell\}$ by its maximal divisible subgroup. 
Also, we denote by $B^{(\ell)}$ the $\ell$-adic completion of $B$, i.e. 
the projective limit
$\underset{^{\longleftarrow}}{{\rm lim}} \, B/\ell^nB$.

For a topological group $B$, we denote by $B^{\wedge}$ the
completion of $B$ with respect to {\em open} subgroups of finite
index (in the discrete case this is the usual profinite completion of
$B$). We set $B^D$ for the group of {\em continuous} homomorphisms
$B\to\Q/\Z$ (in the discrete case these are just all homomorphisms).
We equip $B^D$ with the compact-open topology.
The topological group $B$ is {\em compactly generated} if $B$
contains a compact subset $K$ such that $K$ generates $B$ as a group. 
A continuous morphism $f:\, B\to C$ of topological groups is {\em strict} if 
the image of any open subset of $B$ is an open subset of ${\rm Im} \, f$
for the topology induced by $C$.

\section{Preliminaries on 1-motives} \label{onesect}

Let $S$ be a scheme. Denote by $\fp_S$ the category of {\em fppf} sheaves
of abelian groups over $S$ (when $S=\spec K$ is the spectrum of a field, 
we shall write $\fp_K$ for $\fp_S$). Write  ${\cal C}^b({\fp_S})$ for
the category of bounded complexes of {\em fppf} 
sheaves over $S$ and ${\cal D}^b({\fp_S})$ for the associated derived category.
Recall (e.g. from \cite{ray}) 
that a {\em 1-motive} $M$ over $S$ consists of the following data~:

\begin{itemize}

\item An $S$-group scheme $Y$ which is \'etale locally
isomorphic to $\Z ^r$ for some $r \geq 0$. 

\item A commutative $S$-group scheme $G$ fitting into an exact sequence
of $S$-groups
$$0 \ra T \ra G \stackrel{p}{\ra} A \ra 0$$
where $T$ is an $S$-torus and $A$ an abelian scheme over $S$.

\item An $S$-homomorphism $u : Y \ra G$.

\end{itemize}

The 1-motive $M$ can be viewed as a complex of {\em fppf} $S$-sheaves 
$[Y \stackrel{u}{\ra} G]$, with $Y$ put in degree -1 and $G$ in degree 0,
and also as an object of the derived category ${\cal D}^b(\fp_S)$. It is equipped with a 3-term weight filtration: $W_i(M)=0$ for $i\leq -3$, $W_{-2}(M)=[0\to T]$, $W_{-1}(M)=[0\to G]$ and $W_i(M)=M$ for $i\geq 0$. From this we shall only need the 
1-motive $M/W_{-2}(M)$, i.e. the complex $[Y \stackrel{h}{\ra} A]$,
where $h=p \circ u$.
By \cite{ray}, Proposition~2.3.1, we can identify morphisms 
of 1-motives in ${\cal C}^b({\fp_S})$ and ${\cal D}^b({\fp_S})$.

\medskip

To each 1-motive $M$ one can associate a {\em Cartier dual} $M^{\ast}$ by the following construction which we briefly recall. Denote by $Y^{\ast}$ the group of characters of $T$,
by $A^{\ast}$ the abelian scheme dual to $A$, and by $T^{\ast}$ the 
$S$-torus with character group $Y$. According to the generalised Barsotti-Weil formula (\cite{oort}, \Romannumeral3.18), $A^{\ast}$ represents the functor
$S' \mapsto \ext_{S'}(A,\G)$ on $\fp_S$. Writing $M'=M/W_{-2}M$, one deduces from this that  the functor 
$S' \mapsto \ext_{S'}(M',\G)$ on ${\cal C}^b({\fp_S})$ is representable
by an $S$-group scheme $G^{\ast}$ which is an extension of $A^{\ast}$ by
$T^{\ast}$. One calls the 1-motive $[0\to G^{\ast}]$ the {\em (Cartier) dual of $M'$}. Pulling back the {\em Poincar\'e biextension} (\cite{adt}, p. 395) on $A\times A^{\ast}$ to $A\times G^{\ast}$ one gets a {\em biextension} ${\cal P}'$ of the 1-motives $M'$ and $[0\to G^{\ast}]$ by $\G$, which is a $\G$-torsor over $A\times G^{\ast}$ whose pullback to $Y\times G^{\ast}$ by the natural map $Y\times G^{\ast}\to A\times G^{\ast}$ is trivial (cf. \cite{deligne}, 10.2.1 for this definition).

To treat the general case, consider $M$ as an extension of $M'$ by $T$. Any element of $Y^{\ast}=\Hom_S(T,\G)$ 
then induces by pushout an extension of $M'$ by $\G$, i.e. an element of $G^{\ast}$. Whence a map $Y^{\ast}\to G^{\ast}$ which is in fact a map of $S$-group schemes; we call the associated 1-motive $M^{\ast}$ the {\em (Cartier) dual of $M$}. The pullback $\cal P$ of the biextension ${\cal P}'$ from $A\times G^{\ast}$ to $G\times G^{\ast}$ becomes trivial over $G\times Y^{\ast}$ when pulled back by the map $G\times Y^{\ast}\to G\times G^{\ast}$, hence defines a biextension of $M$ and $M^{\ast}$ by $\G$ (again in the sense of \cite{deligne}).

According to the formula 
$$
{\rm Biext}_S(M, M^{\ast}, \G)\cong\Hom_{{\cal D}^b({\fp_S})}(M\otimes^{\bf L} M^{\ast}, \G[1]) 
$$
of (\cite{deligne}, 10.2.1), the biextension $\cal P$ defines a map
\begin{equation}\label{derpairing}
\Phi_M : M \otimes^{\bf L} M^{\ast} \ra \G[1]
\end{equation}
in ${{\cal D}^b({\fp_S})}$, whence pairings
\begin{equation} \label{pairing}
\hyp^i(S,M) \times \hyp^j(S,M^{\ast}) \ra H^{i+j+1}(S,\G)
\end{equation}
for each $i,j \geq 0$.

(Except when explicitly specified, the cohomology groups in this
paper are relative to the \'etale topology; here
we can work with either the \'etale or the {\em fppf} topology because $G$ 
is smooth over $S$ and $Y$ is \'etale locally constant).
\medskip

\begin{rema}\rm
Two special cases of this construction are classical~:
\begin{itemize}
\item $Y=A=0$, $M=[0 \ra T]$, $M^{\ast}=Y^{\ast}[1]$. Then the pairing 
(\ref{pairing}) is  just the cup-product $H^r(S,T) \times H^{s+1}(S,Y^{\ast})
\ra H^{r+s+1}(S,\G)$. (Similarly for $T=A=0$, $M=Y[1]$, $M^{\ast}=T^{\ast}$.) 

\item $Y=T=0$, $M=[0 \ra A]$, 
$M^{\ast}=[0 \ra A^{\ast}]$. Then (\ref{pairing}) is the well-known pairing in the cohomology of abelian varieties coming from
the generalised Barsotti-Weil formula (compare \cite{adt}, p. 243 and Chapter \Romannumeral3, Appendix C). 
\end{itemize}
\end{rema}

Recall finally (\cite{deligne}, 10.1.5 and 10.1.10) that for any integer $n$ invertible on $S$ and 1-motive $M$ one has an ``$n$-adic realisation'', namely the finite sheaf (or group scheme) defined by
$$
T_{\Z/n\Z}(M)=H^0(M[-1]\otimes^{\bf L} \Z/n\Z),
$$ 
which can be explicitly calculated using the flat resolution $[\Z\mathopto\limits^n \Z]$ of $\Z/n\Z$. The pairing  (\ref{derpairing}) then induces a perfect pairing
\begin{equation}\label{weilpairing}
T_{\Z/n\Z}(M)\otimes T_{\Z/n\Z}(M^{\ast})\to\mu_n
\end{equation}
where $\mu_n$ is the sheaf of $n$-th roots of unity. The classical case is when $M$ is of the form $[0\to A]$ with $A$ an abelian variety, where we find the well-known Weil pairing.

We finish this section by introducing some notation: for each prime number $\ell$ invertible on $S$, we denote by
$T(M)\{ \ell \}$ the direct limit of the $T_{\Z/\ell^n\Z}(M)$ over all $n>0$ and by $T_{\ell}(M)$ their inverse limit. 
The piece of notation $T(M)_{\rm tors}$ stands for the direct sum (taken over all primes $\ell$ invertible on $S$) of the groups $T_{\Z/\ell^n\Z}(M)$.

\section{Local results} \label{deuxsect}

In this section $S$ is the spectrum of a field $K$, complete with respect to a 
discrete valuation and with finite residue field. In particular $K$ 
is a $p$-adic field if ${\rm char } \, K=0$, and is isomorphic to 
the field $\F_q((t))$ for some finite field $\F_q$ if ${\rm char } \, K >0$.
We let $\calo_K$ denote the ring of integers of $K$ and $\F$ its residue field.
\begin{lem} For a 1-motive $M=[Y \ra G]$ over $K$, we have
\begin{itemize}
\item $\hyp^{-1}(K,M)\cong{\rm Ker} \, [H^0(K,Y) \ra H^0(K,G)]$, 
a finitely generated free abelian group; 
\item $\hyp^2(K,M)\cong{\rm Coker} \, [H^2(K,Y) \ra H^2(K,G)]$;
\item $\hyp^i(K,M)=0 \quad \quad i \neq -1,0,1,2$.
\end{itemize}
\end{lem}

\begin{dem}
The field $K$ has strict Galois cohomological dimension 2 (\cite{adt},
\Romannumeral1.1.12).
Since $G$ is smooth, $H^i(K,G)=0$ for any $i >2$; by \cite{adt},
\Romannumeral1.2.1, we also have $H^i(K,Y)=0$ for $i >2$, whence
the last equality. 
For the first two, use moreover  
the distinguished triangle
\begin{equation} \label{triangle}
Y \ra G \ra M \ra Y[1]
\end{equation}
in ${\cal C}^b({\fp_K})$.
\end{dem}

Using the trace isomorphism $H^2(K,\G)\cong\Q/\Z$ of local class field theory, the pairing (\ref{pairing}) of the previous section induces  bilinear pairings
\begin{equation} \label{locpair}
\hyp^i(K,M) \times \hyp^{1-i}(K,M^{\ast}) \ra \Q / \Z
\end{equation}
for all integers $i$ (by the previous lemma, they are trivial for $i\neq -1,0,1,2$).
\bigskip

For $i=-1,1,2$,  we endow the group $\hyp^i(K,M)$ with 
the discrete topology. To topologize $\hyp^0(K,M)$ we proceed as follows. The exact
triangle (\ref{triangle}) yields an exact sequence of abelian groups
\begin{equation} \label{longue}
0 \ra L \ra G(K) \ra \hyp^0(K,M) \ra H^1(K,Y) \ra H^1(K,G)
\end{equation}
where $L:=H^0(K,Y)/\hyp^{-1}(K,M)$ is a discrete abelian group of finite type. 
We equip $I=G(K)/{\rm Im}\,(L)$ with the quotient topology (note that in general 
it is not Hausdorff). The cokernel of the map $G(K)\ra\hyp^0(K,M)$ being finite (as $H^1(K,Y)$ itself is finite by \cite{cogal}, \Romannumeral2.5.8 $iii)$), 
we can define a natural topology on 
$\hyp^0(K,M)$ by taking as a basis of open neighbourhoods of zero the open neighbourhoods of zero in $I$ (this makes $I$ an open subgroup of finite index 
of $\hyp^0(K,M)$).  \medskip

Already in the classical duality theorem for tori over local fields one has to take the profinite completion on $H^0$ in order to obtain a perfect pairing.
However, for the generalisations we have in mind a nuisance arises from the fact that the completion functor is not always left exact, 
even if one works only with discrete lattices and $p$-adic Lie groups.
As a simple example, consider $K=\Q_p$ ($p \geq 3$) 
and the injection $\Z\hookrightarrow\Q_p^{\times}$ given by sending $1$ to $1+p$. Here the induced map on completions $\widehat \Z \ra (\Q_p^{\times})^{\wedge}$ is not injective
(because $\Q_p^{\times} \simeq \Z \times \F_p^{\times} 
\times\Z_p$ and the image of $\Z$ lands in the $\Z_p$-component).\medskip

 Bearing this in mind, for a 1-motive $M=[Y\to G]$ we denote by $\hyp^{-1}_{\wedge}(K,M)$ the kernel of the 
map $H^0(K,Y)^{\wedge} \ra H^0(K,G)^{\wedge}$ coming from $Y \ra G$. There is always an injection $\hyp^{-1}(K,M)^{\wedge}\to \hyp^{-1}_{\wedge}(K,M)$ but it is not an isomorphism in general; the previous example comes from the 1-motive $[\Z\to\G]$.\medskip

However, we shall also encounter a case where the completion functor behaves well.

\begin{lem}\label{complete}
Let $G$ be a semi-abelian variety over the local field $K$, with abelian quotient $A$ and toric part $T$. Then the natural sequence
$$
0\to T(K)^{\wedge}\to G(K)^{\wedge}\to A(K)^{\wedge}\to H^1(K,T)^{\wedge}
$$ is exact. Moreover, $G(K) \hookrightarrow G(K)^{\wedge}$ and 
$(G(K)^{\wedge})^D=G(K)^D$.
\end{lem}

Here in fact we have $A(K)^{\wedge}=A(K)$ (the group $A(K)$ being compact and 
completely disconnected, hence 
profinite) and $H^1(K, T)^{\wedge}=H^1(K, T)$ by finiteness of $H^1(K, T)$ 
(\cite{adt}, \Romannumeral1.2.3).

\dem To begin with, 
the maps between completions are well defined because the maps $T(K) \ra G(K)$, 
$G(K) \ra A(K)$, and $A(K) \ra H^1(K,T)$ are continuous (by \cite{margulis}, \Romannumeral1.2.1.3, 
$T(K)$ is closed 
in $G(K)$ 
and the image of $G(K)$ is open in $A(K)$ by the implicit function theorem). 
The theory of 
Lie groups over a local field 
shows that $G(K)$ is locally compact, completely 
disconnected, and compactly generated; we conclude with the third part of the proposition proven
in the appendix.
\enddem

Now we can state the main result of this section.

\begin{theo} \label{first}
Let $M=[Y \ra G]$ be a 1-motive over the local field $K$. 
The pairing (\ref{locpair}) induces a perfect pairing between
\begin{enumerate}
\item the profinite group $\hyp^{-1}_{\wedge}(K,M)$ and the discrete group 
$\hyp^2(K,M^{\ast})$;
\item the profinite group $\hyp^0(K,M)^{\wedge}$ and the discrete group 
$\hyp^1(K,M^{\ast})$.
\end{enumerate}
\end{theo}
In the special cases $M=[0\to T]$ or $M=[Y\to 0]$ we recover 
Tate-Nakayama duality for tori over $K$ (\cite{cogal}, \Romannumeral2.5.8 and 
\cite{adt}, \Romannumeral1.2.3 for the positive characteristic case) 
and in the case $M=[0\to A]$ we recover Tate's $p$-adic duality theorem for 
abelian varieties and its generalisation to the positive characteristic case
due to Milne (\cite{wc}, \cite{adt}, Cor. \Romannumeral1.3.4, and Theorem 
\Romannumeral3.7.8).

\begin{dem}
For the first statement, set $M':=M/W_{-2} M$. The dual of $M'$ is of the form $[0 \ra G^{\ast}]$, where $G^{\ast}$ is an extension of $A^{\ast}$ by $T^{\ast}$. 
Via the pairing (\ref{locpair}) for $i=-1,0$, we obtain 
a commutative diagram

$$
\begin{CD}
0 @>>> \hyp^{-1}_{\wedge}(K,M') @>>> H^0(K,Y)^{\wedge} @>>> 
H^0(K,A)^{\wedge} \\
&& @VVV @VVV @VVV \\
0 @>>> H^2(K,G^{\ast})^D @>>> H^2(K,T^{\ast})^D @>>> 
H^1(K,A^{\ast})^D 
\end{CD} 
$$
The first line of this diagram is exact by definition, and the second 
one is exact because it is the dual of an exact sequence of 
discrete groups (recall that $H^2(K,A^{\ast})=0$ by \cite{cogal}, \Romannumeral2.5.3, Prop. 16 and \cite{adt}, \Romannumeral3.7.8).
By Tate duality for abelian varieties and Tate-Nakayama duality for tori, the last two vertical 
maps are isomorphisms, hence the same holds for the first one. 

\smallskip

Now using Lemma \ref{complete} we get that the map $H^0(K,Y) \ra 
H^0(K,G)$ induces a map $\hyp^{-1}_{\wedge}(K,M') \ra 
T(K)^{\wedge}$ with kernel $\hyp^{-1}_{\wedge}(K,M)$. 
>From the definition of $M'$
we get a commutative diagram with exact rows
$$
\begin{CD}
0 @>>> \hyp^{-1}_{\wedge}(K,M) @>>> \hyp^{-1}_{\wedge}(K,M') @>>>
H^0(K,T)^{\wedge} \\
&& @VVV @VVV @VVV \\
0 @>>> \hyp^2(K,M^{\ast})^D  @>>> H^2(K,G^{\ast})^D @>>> H^2(K,Y^{\ast})^D  
\end{CD}
$$
whence we conclude as above that the left vertical map is an isomorphism, by the first part and Tate-Nakayama duality. Then $H^2(K, M^*)\cong H^{-1}_{\wedge}(K, M)^D$ follows by dualising, using the isomorphism $H^2(K, M^{\ast})^{DD}\cong H^2(K, M^*)$ for the discrete torsion group $H^2(K, M^{\ast})$.

\medskip

For the second statement, we also begin by working with $M'$. Using the pairings (\ref{locpair}) and Lemma \ref{complete} (applied to $G^{\ast}$), we get a commutative diagram with exact rows:
$$
\begin{CD}
0 @>>>\! T^{\ast}(K)^{\wedge} @>>>\! G^{\ast}(K)^{\wedge} @>>>\! A^{\ast}(K)^{\wedge} @>>>\! H^1(K,T^{\ast})^{\wedge}\\
&& @VVV @VVV @VVV @VVV  \\
0 @>>>\! H^2(K, Y)^D @>>>\! \hyp^1(K,M')^D @>>>\! H^1(K, A)^D @>>>\! H^1(K, Y)^D
\end{CD}
$$
Using the local dualities for $(A, A^{\ast})$ and $(Y, T^{\ast})$,
this implies that
the map $G^{\ast}(K)^{\wedge}\to \hyp^1(K,M')^D$ is an isomorphism.

\smallskip

Now the distinguished triangle $T\to M\to M' \to T[1]$ in 
${\cal C}^b({\fp_K})$
induces the following commutative diagram with exact rows:
$$ 
\begin{CD}
H^0(K,Y^{\ast})^{\wedge} @>>> G^{\ast}(K)^{\wedge} @>>> \hyp^0 (K, M^{\ast})^{\wedge} @>>> H^1(K, Y^{\ast}) @>>> H^1(K, G^{\ast})\\
@VVV @VVV @VVV @VVV @VVV \\
H^2(K, T)^D @>>> \hyp^1(K, M')^D @>>> \hyp^1(K, M)^D @>>> H^1(K, T)^D @>>> \hyp^0(K, M')^D
\end{CD}
$$

Here the exactness of the rows needs some justification. 
The upper row is exact without completing the first three terms. 
Completion in the third term is possible by finiteness of the fourth, 
and completion in the first two terms is possible because the map 
$G^{\ast}(K)\to\hyp^0(K, M^{\ast})$ is open with finite cokernel 
by definition of the topology on the target. In the lower row dualisation 
behaves well because the first four terms are duals of discrete torsion groups. \medskip

By Tate-Nakayama duality for tori and what we have already proven, the first, second and fourth vertical maps are isomorphisms. To derive an isomorphism in the middle it remains to prove the injectivity of the fifth map.  
\smallskip

This in turn follows from the commutative diagram with exact rows (where again we have used the finiteness of $H^1(K, T^{\ast})$ and of $H^1(K, Y)$):

$$
\begin{CD}
A^{\ast}(K)@>>> H^1(K, T^{\ast}) @>>> H^1(K, G^{\ast}) @ >>> H^1(K, A^{\ast})\\
@VVV @VVV @VVV @VVV \\
H^1(K, A)^D @>>> H^1(K, Y)^D @>>> \hyp^0(K, M')^D @>>>  A(K)^D
\end{CD}
$$
Here the first, second and fourth vertical maps are isomorphisms by local duality for tori and abelian varieties. Again,  exactness at the third term of the lower row follows from the definition of the topology on $\hyp^0(K, M')$.
Finally the map $\hyp^0 (K, M^{\ast})^{\wedge} \ra \hyp^1(K, M)^D$ is an 
isomorphism and applying this statement to $M^{\ast}$ instead of $M$, 
we obtain the theorem.
\end{dem}

\bigskip
\bigskip

\begin{rema}\label{charzero}
{\rm 
If $K$ is of characteristic zero, any subgroup of finite index of $T(K)$ 
is open (cf. \cite{adt}, p.32). 
It is easy to see that in this case $\hyp^0(K,M)^{\wedge}$ 
is just the profinite completion of $\hyp^0(K,M)$.
}
\end{rema}

\medskip

\noindent Next we state a version of
Theorem \ref{first} for henselian fields that will be needed for the
global theory.

\begin{theo}\label{hensel}
Let $F$ be the field of fractions of a henselian discrete valuation ring 
$R$ with finite residue field and let $M$ be a 1-motive over $F$. Assume 
that $F$ is of characteristic zero. Then 
the pairing (\ref{locpair}) induces perfect pairings
$$\hyp^{-1}_{\wedge}(F,M) \times \hyp^2(F,M^{\ast}) \ra \Q/\Z$$
$$\hyp^0(F,M)^{\wedge} \times \hyp^1(F,M^{\ast}) \ra \Q/\Z$$
where $\hyp^{-1}_{\wedge}(F,M):={\rm Ker} \, [H^0(F,Y)^{\wedge}
\ra G(F)^{\wedge}]$ and $^{\wedge}$ means profinite completion.
\end{theo}

\begin{remas}\rm ${}$
\begin{enumerate}
\item Denoting by $K$ the completion of $F$, the group $\hyp^0(F,M)$ injects into 
$\hyp^0(K,M)$ by the lemma below, hence it is natural 
to equip $\hyp^0(F,M)$ with the topology induced by $\hyp^0(K,M)$. 
But we shall also show that $\hyp^0(F,M)$ and $\hyp^0(K,M)$ have the same 
profinite completion, hence by Remark~\ref{charzero} the profinite 
completion of $\hyp^0(F,M)$ coincides with its completion with respect
to open subgroups of finite index. Therefore there is no incoherence 
in the notation.
\item In characteristic $p >0$, the analogue of Theorem~\ref{hensel} is not 
clear because of the $p$-part of the groups. Compare \cite{adt}, 
\Romannumeral3.6.13.
\end{enumerate}
\end{remas}

Taking the first remark into account, Theorem~\ref{hensel} immediately results from Theorem \ref{first} via the following lemma.

\begin{lem}\label{green}
Keeping the assumptions of the theorem, denote by $K$ the completion of
$F$. Then the natural map $\hyp^i(F,M) \ra \hyp^i(K,M)$ is an injection for $i=0$ inducing an isomorphism $\hyp^0(F,M)^{\wedge} \ra \hyp^0(K,M)^{\wedge}$ on completions, and an isomorphism 
for $i \geq 1$. 
\end{lem}

\begin{dem}
For any $n>0$, the canonical map $G(F)/n\to G(K)/n$ is surjective, for $G(F)$
is dense in $G(K)$ by Greenberg's approximation theorem \cite{greenberg},
and $nG(K)\subset G(K)$ is an open subgroup. But this map is also injective,
for any  point 
$P\in G(K)$ with $nP\in G(F)$ is locally given by coor\-dinates algebraic 
over $F$, but $F$ is algebraically closed in $K$
(apply e.g. \cite{nagata}, Theorem 4.11.11 and note that $F$ is of
characteristic 0), hence $P\in G(F)$. 
Since $Y$ is locally constant in the \'etale topology over $\spec F$,
we have $H^i(F,Y)=H^i(K,Y)$ for each $i \geq 0$. The case $i=0$ of the lemma follows from these facts by d\'evissage. 

To treat the cases $i>0$, recall first that multiplication 
by $n$ on $G$ is surjective in the \'etale topology. Therefore
$$H^i(F,G)[n]=\cok [H^{i-1}(F,G)/n \ra H^{i}(F,G[n])]$$
for $i \geq 1$, and similarly for $H^i(K,G)$. 
Moreover, $H^i(F,G[n])=H^i(K,G[n])$ because 
$G[n]$ is locally constant in the \'etale topology (note that $F$ and $K$
have the same absolute Galois group). Starting from the isomorphism $G(F)/n \cong G(K)/n$ already proven, we thus obtain
isomorphisms of torsion abelian groups
$H^i(F,G) \simeq H^i(K,G)$ for any $i \geq 1$ by induction on $i$, which together with the similar isomorphisms for $Y$ mentioned above yield the statement by d\'evissage. 
\end{dem}

We shall also need the following slightly finer statement.

\begin{prop}\label{greenbis}
Keeping the notations above, equip $\hyp^0(F,M^{\ast})$ 
with the topology induced by $\hyp^0(K,M^{\ast})$.
Then $\hyp^0(F,M^{\ast})$ injects 
into $\hyp^0(F,M^{\ast})^{\wedge}$ and they have the same continuous dual. Moreover, the pairing (\ref{locpair}) yields an isomorphism
$$\hyp^1(F,M) \cong \hyp^0(F,M^{\ast})^D.$$
\end{prop}

\begin{dem} For the first statement we use the exact sequence 
$$H^0(F,Y^{\ast}) \ra G^{\ast}(F) \ra \hyp^0 (F, M^{\ast}) \ra 
H^1(F, Y^{\ast}) \ra H^1(F, G^{\ast})$$
and similarly for $K$ instead of $F$. By definition of the topology 
on $\hyp^0 (F, M^{\ast})$ and $\hyp^0 (K, M^{\ast})$, the duals of 
these sequences remain exact thanks to the result proven in the Appendix. The group $G^{\ast}(F)$ injects into $G(K)^{\wedge}$ by Lemma \ref{green} and Lemma \ref{complete}; moreover,  
since $G^{\ast}(F)$ is a dense subgroup of $G^{\ast}(K)$, they have the same dual. Similarly, we have an injection $H^0(F, Y^{\ast})\hookrightarrow H^0(K, Y^{\ast})^{\wedge}$ using that $H^0(K, Y^{\ast})$ is of finite type, and this map induces an isomorphism on duals. Thus the first statement follows from the exact sequence using the finiteness of $H^1(F, Y^{\ast})\hookrightarrow H^1(K, Y^{\ast})$ and Lemma \ref{green}. The second statement now follows from Theorem \ref{first}, again using Lemma \ref{green}.
\end{dem}

For the global case, we shall also need a statement for the real case.
Consider a 1-motive $M_{\RR}$ over the spectrum of 
the field $\RR$ of real numbers. As in the classical cases, the duality
results in the previous section extend in a straightforward fashion to
this situation, provided that we replace usual Galois cohomology groups
by Tate modified groups. Denote by $\Gamma_{\RR}=\gal(\C/\RR) \simeq 
\Z/2$ the Galois group of $\RR$.
Let ${\cal F}^{\bullet}$ be a bounded complex 
of $\RR$-groups.  For each $i \in \Z$, the 
modified hypercohomology groups $\thyp^i(\RR,{\cal F}^{\bullet})$ are defined 
in the usual way: for each term ${\cal F}^i$ 
of ${\cal F}$, we take the standard Tate complex associated to the 
$\Gamma_{\RR}$-module ${\cal F}^i(\C)$ (cf. \cite{adt}, pp. 2--3); then
we obtain Tate hypercohomology groups via the complex associated to the 
arising double complex. From the corresponding well-known results 
in Galois cohomology, it is easy to see that 
$\thyp^i(\RR,{\cal F}^{\bullet})=\hyp^i(\RR,{\cal F}^{\bullet})$ for 
$i \geq 1$ if ${\cal F}^{\bullet}$ is concentrated in nonpositive degrees, and that $\thyp^i(\RR,{\cal F}^{\bullet})$ is isomorphic to 
$\thyp^{i+2}(\RR,{\cal F}^{\bullet})$ for any $i \in \Z$. 
Recall also that the Brauer group $\br {\RR}$ is isomorphic to 
$\Z/2\Z \subset \Q/\Z$ via the local invariant.

\smallskip

Now we have the following analogue of Theorem~\ref{first}:

\begin{prop} \label{realprop}
Let $M_{\RR}=[Y_{\RR} \ra G_{\RR}]$ be a 1-motive over $\RR$. Then the 
cup-product pairing 
induces a perfect pairing of finite 2-torsion groups
$$\thyp^0(\RR,M_{\RR}) \times \thyp^1(\RR,M_{\RR}^{\ast}) \ra \Z/2\Z$$
\end{prop}

\dem Let $T_{\RR}$ (resp. $A_{\RR}$) be the torus (resp. the abelian 
variety) corresponding to $G_{\RR}$. 
In the special cases $M_{\RR}=T_{\RR}$, $M_{\RR}=A_{\RR}$, 
$M_{\RR}=Y_{\RR}[1]$, the result is known (\cite{adt}, \Romannumeral1.2.13 
and \Romannumeral1.3.7). Now the proof by devissage consists exactly of 
the same steps as in Theorem~\ref{first}, except that we don't have to
take any profinite completions, all occurring groups being finite.
\enddem

In Section~5 
we shall need the fact that when $M$ is a 1-motive over a local field $K$
which extends to a 1-motive over
$\spec \calo_K$, the unramified parts of the cohomology are exact annihilators
of each other in the local duality pairing for $i=1$ (see \cite{cogal}, \Romannumeral2.5.5, \cite{neukirch}, Theorem 7.2.15 and 
\cite{adt}, \Romannumeral3.1.4 for analogues for finite modules). More
precisely, let ${\cal M}=[{\cal Y} \ra {\cal G}]$ 
be a 1-motive over $\spec \calo_K$ and $M=[Y \ra G]$ the 
restriction of ${\cal M}$ to $\spec K$. Denote by $\hyp^0_{\nr}(K,M)$
and $\hyp^1_{\nr}(K,M^{\ast})$ the respective images of the maps
$\hyp^0(\calo_K,{\cal M})\to\hyp^0(K,M)$ and $\hyp^1(\calo_K,{\cal M}^{\ast})\to\hyp^1(K,M^{\ast})$. 
To make the notation simpler, we still let $\hyp^0_{\nr}(K,M)^{\wedge}$
denote the image of $\hyp^0_{\nr}(K,M)^{\wedge}$ in $\hyp^0(K,M)^{\wedge}$. (We
work with complete fields since this is what will be needed later; the
henselian case is similar in mixed characteristic.)

\begin{theo} \label{unram}
In the above situation, $\hyp^0_{\nr}(K,M)^{\wedge}$ and
$\hyp^1_{\nr}(K,M^{\ast})$ are the exact annihilators 
of each other in the pairing $$\hyp^0(K,M)^{\wedge} \times 
\hyp^1(K,M^{\ast}) \ra \Q/\Z$$ induced by (\ref{locpair}).
\end{theo}

\dem  The restriction of the local pairing to $\hyp^0_{\nr}(K,M) 
\times \hyp^1_{\nr}(K,M^{\ast})$ is zero because 
$H^2(\calo_K,\G)\cong H^2({\bf F}, \G)=0$. Thus it is sufficient to show that the maps 
$$\hyp^0(K,M)^{\wedge}/\hyp^0_{\nr}(K,M)^{\wedge} 
\ra \hyp^1_{\nr}(K,M^{\ast})^D,$$
$$\hyp^1(K,M^{\ast})/\hyp^1_{\nr}(K,M^{\ast}) \ra \hyp^0_{\nr}(K,M)^D$$
are injective, where we have equipped 
$\hyp^1_{\nr}(K,M^{\ast})$ with the discrete topology
and $\hyp^0_{\nr}(K,M)$ with the topology induced by that on
$\hyp^0(K,M)$.  

\noindent Denote by ${\cal T}$ (resp. $T$) the torus
and by ${\cal A}$ (resp. $A$) the 
abelian scheme (resp. abelian variety) corresponding 
to ${\cal G}$ (resp. $G$). We need the following lemma presumably well known to
the experts.

\begin{lem} \label{fppf} In the Tate-Nakayama pairing $$H^2(K,Y) \times H^0(K,T^{\ast})\ra \Q/\Z,$$ 
the exact annihilator of $H^0_{\nr}(K,T^{\ast})$ is $H^2_{\nr}(K,Y)$.
\end{lem}

\dem Let $n >0$. We work in flat cohomology. The exact sequence of {\em fppf} sheaves
$$0 \ra {\cal T}[n] \ra {\cal T} \ra {\cal T} \ra 0$$
and the cup-product pairings induce a commutative diagram with exact rows:

$$
\begin{CD}
&& H^1_{fppf}(\calo_K,{\cal Y}/n{\cal Y}) @>>> 
H^2_{fppf}(\calo_K,{\cal Y})[n] \\
&& @VVV @VVV \\
&& H^1_{fppf}(K,Y/nY) @>>> H^2_{fppf}(K,Y)[n] @>>>0 \\
&& @VVV @VVV \\
0 @>>> H^1_{fppf}(\calo_K, {\cal T}^{\ast}[n])^D @>>> 
H^0(\calo_K,{\cal T}^{\ast})^D 
\end{CD}
$$
where all groups 
are given the discrete topology. The zero at lower left comes from the
vanishing $H^1_{fppf}(\calo_K,{\cal T}^{\ast} )\cong H^1_{fppf}({\bf F},
{\widetilde T}^{\ast})=0$ (where ${\widetilde T}^{\ast}$ stands for the 
special fibre of ${\cal T}^{\ast}$),
which is a consequence of Lang's theorem (\cite{lang}, Theorem 2 
and \cite{milne}, \Romannumeral3.3.11).

\smallskip

\noindent Now by \cite{adt}, \Romannumeral3.1.4. and \Romannumeral3.7.2., 
the left column is exact.
Therefore the right column is exact as well. To see that this implies the
statement it remains to note that since $\cal Y$ is a smooth group scheme over $\spec\calo_K$, its \'etale and flat
cohomology groups are the same and moreover they are all torsion in positive degrees. 
\enddem

\smallskip We resume the proof of Theorem~\ref{unram}. 
 The weight filtration on $M$, the cup-product pairings and the inclusion
$\calo_K\subset K$  induce a commutative diagram with exact rows (here the groups 
in the lower row are given the discrete topology)
\begin{equation}\label{funny}
\begin{CD}
&&&& \hyp^1(\calo_K,{\cal M}) @>>> \hyp^2(\calo_K,{\cal Y}) @>>> 0 \\
&&&& @VVV @VVV \\
&& H^1(K,G) @>>> \hyp^1(K,M)  @>>> H^2(K,Y) \\
&& @VVV @VVV @VVV \\
0 @>>> \hyp^0(\calo_K, [{\cal Y}^{\ast} \ra {\cal A}^{\ast}])^D
@>>> \hyp^0(\calo_K,{\cal M}^{\ast})^D @>>> H^0(\calo_K, {\cal T}^{\ast})^D
\end{CD}
\end{equation}
where the two zeros come from the vanishing of the groups 
$H^1(\calo_K,{\cal T}^{\ast})$ and
$H^2(\calo_K,{\cal G})=H^2(\F,{\widetilde G})$; the second vanishing follows from the fact that ${\widetilde G}(\bar \F)$ is torsion and $\F$ is of cohomological dimension 1. 
\smallskip

Next, observe that the map 
$\hyp^0(\calo_K,[{\cal Y}^{\ast} \ra {\cal A}^{\ast}]) \ra 
\hyp^0(K,[Y^{\ast} \ra A^{\ast}])$ is an isomorphism. Indeed, by d\'evissage
this reduces to showing that the natural maps $H^0(\calo_K,{\cal A})\to
H^0(K,A)$ and $H^1(\calo_K,{\cal Y})\to
H^1(K,Y)$ are isomorphisms. The first isomorphism follows from the  
properness of the abelian scheme ${\cal A}$. For the second, denote by $\calo_K^{nr}$ the strict henselisation of
$\calo_K$ and by $K^{nr}$ its fraction field. Then the Hochschild-Serre spectral sequence induces a commutative diagram with exact rows:
$$
\begin{CD}
0 @>>> H^1({\bf F}, H^0(\calo_K^{nr}, {\cal Y}))@>>> H^1(\calo_K, {\cal Y})@>>> H^0({\bf F},
H^1(\calo_K^{nr}, {\cal Y}))\\
&& @VVV  @VVV  @VVV \\
0 @>>> H^1({\bf F}, H^0(K^{nr}, Y))@>>> H^1(K, Y)@>>> 
H^0({\bf F}, H^1(K^{nr}, Y))
\end{CD}
$$
Here the group at top right vanishes because $\calo_K^{nr}$ is acyclic for \'etale cohomology. Also, since $\calo_K^{nr}$ is simply connected for the \'etale topology, the sheaf $\cal Y$ is isomorphic to a (torsion free) constant sheaf $\Z^r$, whence the vanishing of the group at bottom right. For the same reason, both groups on the left are isomorphic to $H^1(\gal({\bf F}'/{\bf F}), H^0(\calo_K^{nr}, {\cal Y}))$, where $\bf F'$ is a finite extension trivialising the action of $\gal(\bar F|F)$ on $H^0(\calo_K^{nr}, {\cal Y})$, whence the claim.

This being said, we conclude from
Proposition \ref{greenbis} (which also holds in positive characteristic over
complete fields) that the left vertical map in diagram (\ref{funny}) is injective.
On the other hand the right column is exact by Lemma~\ref{fppf}.
Hence the middle column, i.e. the sequence
$$\hyp^1(\calo_K,M) \ra \hyp^1(K,M) \ra \hyp^0(\calo_K, {\cal M}^{\ast})^D$$
is exact. Since the map $\hyp^1(K,M) \ra \hyp^0(\calo_K, {\cal M}^{\ast})^D$
factors through the map $\hyp^1(K,M) \ra 
(\hyp^0_{\nr}(K,M^{\ast})^{\wedge})^D$, the sequence
$$ \hyp^1(\calo_K,M) \ra \hyp^1(K,M) \ra (\hyp^0_{\nr}(K,M^{\ast})^{\wedge})^D$$
(of which we knew before that it is a complex) 
is exact as well. Dualising this exact sequence of discrete groups, we obtain from
Theorem~\ref{first} that the sequence 
$$\hyp^0_{\nr}(K,M^{\ast})^{\wedge} \ra \hyp^0(K,M)^{\wedge} \ra 
\hyp^1(\calo_K,M)^D$$
is exact and the theorem is proven.
\enddem

\bigskip 

\section{Global results~: \'etale cohomology} \label{threesect}

Let $k$ be a number field with ring of integers $\calo_k$. 
Denote by $\Omega_k$ the set of places of $k$ and 
$\Omega^{\infty}_k\subset\Omega_k$ the subset of real places. 
Let $k_v$ be the completion of $k$ at $v$ if $v$ is archimedean, 
and the field of fractions of the henselisation of the local ring 
of $\spec \calo_k$ at $v$ if $v$ is finite. In the latter case, 
the piece of notation 
${\hat k}_v$ stands for the completion of $k$ at $v$. We denote by $U$ an open
subscheme of $\spec \calo_k $ and by $\Sigma_f$ the set of finite places coming from closed points outside $U$.

In this section and the next, 
every abelian group is equipped with the discrete topology; in
particular $B^{\wedge}$ denotes the profinite completion of $B$ and 
$B^D:=\Hom(B,\Q/\Z)$ (even if $B$ has a natural nondiscrete topology).

\smallskip

We shall need the notion of ``cohomology groups with compact support'' $\hyp^i_c(U, {\cal F}^{\bullet})$ for cohomologically bounded complexes of abelian 
sheaves ${\cal F}^{\bullet}$ on $U$. For $k$ totally imaginary, these satisfy 
$\hyp^i_c(U, {\cal F}^{\bullet})=\hyp^i(\spec \calo_k , 
j_!{\cal F}^{\bullet})$, 
where $j:\, U\to\spec\calo_k$ is the inclusion map. In the general case this equality holds up to a finite 2-group. More precisely, there exists a long exact sequence (infinite in both directions)
{\small
\begin{equation}\label{longex}
\dots\to\hyp^i_c(U, {\cal F}^{\bullet})\to \hyp^i(U, {\cal F}^{\bullet})
\to\bigoplus_{v\in\Sigma_f}\hyp^i(\hat k_v, {\cal F}^{\bullet})
\oplus\bigoplus_{v\in\Omega^{\infty}_k}
\thyp^i(k_v, {\cal F}^{\bullet})\to \hyp^{i+1}_c(U, {\cal F}^{\bullet})\to\dots 
\end{equation}
}
\noindent where for $v\in\Omega_k^{\infty}$ the notation 
$\thyp^i(k_v, {\cal F}^{\bullet})$ stands for Tate (modified) cohomology groups of the group $\gal(\bar k_v/k_v)\cong \Z/2\Z$ and where we have abused notation in denoting the pullbacks of ${\cal F}^{\bullet}$ under the maps $\spec  k_v \to\spec \calo_k $ by the same symbol.

In the literature, two constructions for the groups $\hyp^i_c(U, {\cal F}^{\bullet})$ have been proposed, by Kato \cite{kato} and by Milne \cite{adt}, respectively. The two definitions are equivalent (though we could not find an appropriate reference for this fact). We shall use Kato's construction which we find more natural and which we now copy from \cite{kato} for the convenience of the 
reader. 
 
First, for an abelian sheaf $\cal G$ on the big \'etale site of $\spec \Z $, one defines a complex $\widehat{\cal G}^{\bullet}$ as follows. Denote by $a:\,\spec \C \to\spec \Z $ the canonical morphism and by $\sigma:\, a_{\ast}a^{\ast}{\cal G}\to a_{\ast}a^{\ast}{\cal G}$ the canonical action of the complex conjugation viewed as an element of $\gal(\C/\R)$. Now put $\widehat{\cal G}^0={\cal G}\oplus a_{\ast}a^{\ast}{\cal G}$ and $\widehat{\cal G}^i=a_{\ast}a^{\ast}{\cal G}$ for $i\in\Z\setminus\{0\}$. One defines the differentials $d^i$ of the complex $\widehat{\cal G}^{\bullet}$ as follows: $$d^{-1}(x)=(0, (\sigma-\id)(x));\qquad d^0(x,y)=b(x)+(\sigma+\id)(y),$$ where $b:\, {\cal G}\to a_{\ast}a^{\ast}{\cal G}$ is the adjunction map; otherwise, set $d^i=\sigma+\id$ for $i$ even and $d^i=\sigma-\id$ for $i$ odd. This definition extends to bounded complexes ${\cal G}^{\bullet}$ on the big \'etale site of $\spec \Z $ in the usual way: construct the complex $\widehat{\cal G}^i$ for each term ${\cal G}^i$ of the complex and then take the complex associated to the arising double complex. Finally for $U$ and ${\cal F}^{\bullet}$ as above, one sets
$$
\hyp^i_c(U, {\cal F}^{\bullet}):=\hyp^i(\spec \Z , \widehat{{\bf R}f_!{\cal F}^{\bullet}}),
$$
where $f: U\to\spec \Z $ is the canonical morphism. 

>From this definition one infers that for an open immersion $j_V:\, V\nobreak\to\nobreak U$ and a complex ${\cal F}_V^{\bullet}$ of sheaves on $V$ one has $\hyp^i_c(V, {\cal F}_V^{\bullet})\cong\hyp^i_c(U, j_{V!}{\cal F}_V^{\bullet})$; therefore, setting ${\cal F}_V^{\bullet}=j_V^{\ast}{\cal F}^{\bullet}$ one obtains a canonical map $$\hyp^i_c(V, {\cal F}_V^{\bullet})\to \hyp^i_c(U, {\cal F}^{\bullet})$$ coming from the morphism of complexes $j_{V!}j_V^{\ast}{\cal F}^{\bullet}\to{\cal F}^{\bullet}$. This covariant functoriality for open immersions will be crucial for the arguments in the next section.

Finally, we remark that for cohomologically bounded complexes ${\cal F}^{\bullet}$
and ${\cal G}^{\bullet}$ of \'etale sheaves on $U$, one has a cup-product pairing
\begin{equation}\label{compair}
\hyp^i(U, {\cal F}^{\bullet})\times \hyp^j_c(U, {\cal G}^{\bullet})\to \hyp^{i+j+1}_c(U, {\cal F}^{\bullet}\otimes^{\bf L}{\cal G}^{\bullet}).
\end{equation}

Indeed, for $f$ as above (which is quasi-finite by definition), one knows from general theorems of \'etale cohomology that the complexes ${\bf R}f_*{\cal F}^{\bullet}$ and ${\bf R}f_!{\cal F}^{\bullet}$ are cohomologically bounded (and similarly for ${\cal G}^{\bullet}$, ${\cal F}^{\bullet}\otimes^{\bf L}{\cal G}^{\bullet}$), and that there exists a canonical pairing  ${\bf R}f_*{\cal F}^{\bullet}\otimes^{\bf L}{\bf R}f_!{\cal G}^{\bullet}\to{\bf R}f_!({\cal F}^{\bullet}\otimes^{\bf L}{\cal G}^{\bullet})$. One then uses the simple remark that any derived pairing ${\cal A}^{\bullet}\otimes^{\bf L}{\cal B}^{\bullet}\to{\cal C}^{\bullet}$ of cohomologically bounded complexes of \'etale sheaves on $\spec \Z$ induces a pairing ${\cal A}^{\bullet}\otimes^{\bf L}\widehat{{\cal B}^{\bullet}}\to\widehat{{\cal C}^{\bullet}}$.
   
\begin{rema}\rm
In \cite{zink}, Zink defines modified cohomology groups $\widehat H^i(U, {\cal
F})$ which take the real places into account and satisfy a localisation sequence for
cohomology. He applies this to prove the Artin-Verdier duality theorem for
finite sheaves in the case $U=\spec\calo_k$ (where cohomology and compact support
cohomology coincide). For general $U$, however, one needs the groups $H^i_c(U, {\cal
F})$ which are the same as the groups ${\widehat H}^i(\spec\calo_k, j_!{\cal F})$ in his notation.  
\end{rema}

This being said, we return to 1-motives.

\begin{lem}\label{sorites} Let $M$ be a 1-motive over $U$.
\begin{enumerate}
\item The groups $\hyp^i(U,M)$ are torsion for $i\geq 1$ and so are the groups $\hyp^i_c(U,M)$  for $i\geq 2$. 
\item For any $\ell$ invertible on $U$, 
the groups $\hyp^i(U,M)\{\ell\}$ ($i \geq 1$)
are of finite cotype. Same assertion for the groups $\hyp^i_c(U,M)\{\ell\}$ ($i \geq 2$).
\item The group $\hyp^0(U,M)$ is of finite type.
\end{enumerate}
\end{lem}

\dem For the first part of $(1)$ note that, with the notation
of Section~\ref{onesect}, the group $H^i(U,A)$ is torsion 
(\cite{adt}, \Romannumeral2.5.1), and so are $H^i(U,T)$ and $H^{i+1}(U,Y)$
by \cite{adt}, \Romannumeral2.2.9. The second part then follows from exact sequence (\ref{longex}) and the local facts. 

By this last argument, for $(2)$ it is again enough to prove the first statement. To do so, observe that
$H^i(U,A)\{l\}$ is of finite cotype (\cite{adt}, \Romannumeral2.5.2)
and for each positive integer $n$, there are surjective
maps 
$$H^i(U,Y/l^n Y) \ra H^{i+1}(U,Y)[l^n], \quad H^i(U,T[l^n]) \ra H^i(U,T)[l^n]$$
whose sources are  finite by
\cite{adt}, \Romannumeral2.3.1. 

To prove $(3)$, one first uses for each $n >0$ 
the surjective map from the finite group $H^0(U,Y/nY)$ onto
$H^1(U,Y)[n]$. 
It shows that the finiteness of $H^1(U,Y)$ follows if we 
know that  $H^1(U,Y)=H^1(U,Y)[n]$ for some $n$. Using a standard
restriction-corestriction argument, this follows from the fact 
that $H^1(V,\Z)=0$ for any 
normal integral scheme $V$ (\cite{sga}, \Romannumeral9.3.6 (ii)). 
It remains to note
that the groups $H^0(U,A)=H^0(k,A)$ and $H^0(U,T)$ are of finite type, by the Mordell-Weil Theorem
and Dirichlet's Unit Theorem, respectively 
(for the latter observe that 
$H^0(U,T)$ injects into $H^0(V,T)\cong (H^0(V, \G))^r$, where $V/U$ is an 
\'etale covering trivialising $T$).
\enddem

\begin{rema}\label{soritebis}\rm
The structure of the group $\hyp^1_c(U, M)$ is a bit more complicated: it is an extension of a torsion group (whose $\ell$-part is of finite cotype for all $\ell$ invertible on $U$) by a quotient of a profinite group.
\end{rema}

Now, as explained on p. 159 of \cite{kato}, combining a piece of the long exact sequence (\ref{longex}) for ${\cal F}^{\bullet}=\G$ with the main results of global class field theory yields a canonical (trace) isomorphism $$H^3_c(U, \G)\cong \Q/\Z.$$ Also, we have a natural compact support version
$$
\hyp^i(U, M)\times \hyp^j_c(U, M^{\ast})\to \hyp^{i+j+1}_c(U, \G)
$$
of the pairing (\ref{pairing}), constructed using the pairing (\ref{compair}) above. Combining the two, we get canonical pairings
$$
\hyp^i(U,M) \times \hyp^{2-i}_c(U,M^{\ast}) \ra \Q/\Z
$$
defined for $-1 \leq i \leq 3$. For any prime number $\ell$ invertible on $U$, restricting to $\ell$-primary torsion and modding out by divisible elements (recall the notations from the beginning of the paper) induces pairings
\begin{equation}\label{ellpairing}
\overline{\hyp^i(U,M)}\{\ell\} \times \overline{\hyp^{2-i}_c(U,M^{\ast})}\{\ell\} \ra \Q/\Z.
\end{equation}

\begin{theo} \label{compactsup}
For any 1-motive $M$ and any $\ell$ invertible on $U$, the pairing (\ref{ellpairing}) is non-degenerate for $0 \leq i \leq 2$.
\end{theo}

Note that the two groups occurring in the pairing (\ref{ellpairing})
are finite by Lemma \ref{sorites} $(2)$.

\begin{dem}
This is basically the argument of (\cite{adt}, \Romannumeral2.5.2 (b)). Let $n$ be a power of $\ell$. Tensoring the exact sequence 
$$
0\to\Z\to\Z\to\Z/n\Z\to 0
$$
by $M$ in the derived sense and passing to \'etale cohomology over $U$   induces exact sequences
$$
0\to \hyp^{i-1}(U, M)\otimes\Z/n\Z\to \hyp^{i-1}(U, M\otimes^{\bf L}\Z/n\Z)\to\hyp^i(U, M)[n]\to 0
$$
Now $M\otimes^{\bf L}\Z/n\Z$ viewed as a complex of \'etale sheaves has trivial cohomology in degrees other than $-1$; indeed, with the notation of 
Section 1, the group $Y$ is torsion free and multiplication by $n$ on the group scheme $G$ is surjective in the \'etale topology. Therefore, using the notation of Section 1, we may rewrite the previous sequence as
\begin{equation}\label{nmultiply}
0\to \hyp^{i-1}(U, M)\otimes\Z/n\Z\to H^i(U, T_{\Z/n\Z}(M))\to\hyp^i(U, M)[n]\to 0.
\end{equation}
Write $T(M)\{\ell\}$ for the direct limit of the groups $T_{\Z/n\Z}(M)$ as 
$n$ runs through powers of $\ell$ and $T_{\ell}(M)$ for their inverse limit. 
For each $r \geq 0$, $H^r(U,T_{\ell}(M))$ stands for 
the inverse limit of $H^r(U,T_{\Z/n\Z}(M))$ ($n$ running through powers 
of $l$), and similarly 
for compact support cohomology.
Passing to the direct limit in the above sequence then induces an isomorphism
$$
\overline{H^i(U, T(M)\{\ell\})}\cong\overline{\hyp^i(U, M)}\{\ell\}.
$$
Now by Artin-Verdier duality for finite sheaves 
(\cite{zink}, \cite{adt}, \Romannumeral2.3; see also \cite{mazur} in the totally imaginary case) the first group here is isomorphic 
(via a pairing induced by (\ref{ellpairing})) to the dual of the group 
$H^{3-i}_c(U, T_{\ell}(M^{\ast}))\{\ell\}$.

Working with the analogue of exact 
sequence (\ref{nmultiply}) for compact support cohomology and passing to the 
inverse limit over $n$ using the finiteness of the groups 
$H^{3-i}_c(U, T_{\Z/n\Z}(M^{\ast}))$, we get isomorphisms 
$$
\hyp^{2-i}_c(U, M^{\ast})^{(\ell)}\{\ell\}
\cong
H^{3-i}_c(U, T_{\ell}(M^{\ast}))\{\ell\}
$$ 
using the torsion freeness of the $\ell$-adic Tate module of the group 
$\hyp^{3-i}_c(U, M^{\ast})$. Finally, we have 
$\overline{\hyp^{2-i}_c(U,M^{\ast})}\{\ell\}  \cong \hyp^{2-i}_c(U,M^{\ast})^
{(\ell)}\{\ell\}$
by the results in Lemma \ref {sorites} and Remark \ref{soritebis}.
\end{dem}

\bigskip

>From now on, we shall make the convention (to ease notation) 
that for any archimedean place $v$ and each $i \in \Z$, $\hyp^i(k_v,M)$
means the {\em modified}
group $\thyp^i(k_v,M)$ (In particular it is zero if $v$ is 
complex, and it is a finite 2-torsion group if $v$ is real).

Following \cite{adt}, \Romannumeral2.5, we define for $i \geq 0$
$$D^i(U,M)={\rm Ker}\, [\hyp^i(U,M) \ra \prod_{v \in \Sigma} \hyp^i(k_v,M)]$$
where the finite subset $\Sigma=\Sigma_f \cup \Omega_{\infty} 
\subset \Omega_k$ consists of the
real places and the primes of $\calo_k$ which do not correspond to a
closed point of $U$. For $i\geq 0$ we also have 
$$D^i(U,M)={\rm Im} \, [\hyp^i_c(U,M) \ra \hyp^i(U,M)]$$
by the definition of compact support cohomology. By Lemma \ref{sorites}, $D^1(U,M)$ is a torsion group and
$D^1(U,M)\{\ell\}$ is of finite cotype. Now
Theorem~\ref{compactsup} has the following
consequence:

\begin{cor} \label{shacor}
Under the notation and assumptions of Theorem~\ref{compactsup}, for $i=0,1$ there is
a pairing
\begin{equation} \label{interm}
D^i(U,M)\{\ell\} \times D^{2-i}(U,M^{\ast})\{\ell\} \ra \Q/\Z
\end{equation}
whose left and right kernels are respectively the divisible subgroups
of the two groups.
\end{cor}

\begin{remas}\label{zerorem}\rm ${}$
\begin{enumerate}
\item As $D^0(U,M)$ is of finite type, the divisible part of 
$D^0(U,M)\{ \ell \}$ is zero.
\item In the case $M=[0\to T]$, the corollary gives back (part of) \cite{adt}, Corollary \Romannumeral2.4.7, itself based on the main result of Deninger \cite{deninger} (which we do not use here).
\end{enumerate}
\end{remas}

\dem The proof works along the same pattern as that used in the case $i=1$, $M=[0\to A]$ treated in \cite{adt}, Corollary 
\Romannumeral2.5.3. As the proof given in {\em loc. cit.} contains some inaccuracies, we give a
detailed argument. 
We use the following commutative diagram whose exact rows come from the definition of
the groups $D^i$ and whose vertical maps are induced by the pairings (\ref{ellpairing}) and (\ref{locpair}). 

$$
\begin{CD}
0 @>>> D^i(U,M)\{\ell\} @>>> \hyp^i(U,M)\{\ell\} @>>> \bigoplus_{v \in \Si} 
\hyp^i(k_v,M) \\
&&&& @VVV @VVV \\
0 @>>> D^{2-i}(U,M^{\ast})^D @>>> 
\hyp^{2-i}_c(U,M^{\ast})^D @>>>
\bigoplus_{v \in \Si} \hyp^{1-i}(k_v,M^{\ast})^D
\end{CD}
$$

\noindent The diagram defines a 
map $D^i(U,M)\{\ell\} \ra D^{2-i}(U,M^{\ast})^D$.
The right vertical map is injective for $i=0$ by Theorem \ref{first}, Remark
\ref{charzero}, the first statement in Proposition~\ref{greenbis}, and 
Proposition~\ref{realprop}.
For $i=1$ it is injective by Proposition~\ref{realprop} and
the second statement in Proposition~\ref{greenbis} 
(indeed for $v$ finite,
$\hyp^0(k_v,M^{\ast})$ is now equipped with the discrete topology, 
which is finer than the topology defined in Proposition~\ref{greenbis}).
The middle vertical map has divisible kernel
by Theorem~\ref{compactsup}. Hence the map $D^i(U,M)\{\ell\} \ra
D^{2-i}(U,M^{\ast})^D$ has divisible kernel, and similarly for the induced 
map $D^i(U,M)\{\ell\} \ra D^{2-i}(U,M^{\ast})\{\ell\}^D$ (recall that 
$D^1(U,M)$ and $D^2(U,M)$ are torsion).
Now exchange the roles of $M$ and $M^{\ast}$.
\enddem

Let $A_k$ denote the generic fibre of $A$ 
and 
$$\Sha^1(A_k):={\rm Ker} \, [H^1(k,A_k) \ra \prod_{v \in \Omega_k}
H^1({\hat k}_v,A_k)]$$ its Tate-Shafarevich group.
A special case of Corollary~\ref{shacor} is

\begin{prop} \label{partsha}
Let $M$ be a 1-motive over $U$ and $\ell$ a prime number invertible on $U$.
Assume that $\Sha^1(A_k)\{\ell\}$ and $\Sha^1(A^{\ast}_k)\{\ell\}$ are finite. 
Then the pairing 
$$D^1(U,M)\{\ell\} \times D^1(U,M^{\ast})\{\ell\} \ra \Q/\Z$$
of Corollary~\ref{shacor} is a perfect pairing of finite groups.
\end{prop}
\dem Using Corollary~\ref{shacor}, it is sufficient to prove that 
$D^1(U,M)\{\ell\}$ is finite. We have a commutative diagram with exact rows:
$$
\begin{CD}
\newenvironment{disarray}
{\everymath{\displaystyle\everymath{}}\array}
{\endarray}
H^1(U,Y)@>>> H^1(U, G) @>>> \hyp^1(U, M) 
@>>> H^2(U, Y) \\
@VVV @VVV @VVV @VVV \\
\bigoplus_{v \in \Si} H^1(k_v,Y) @>>> \bigoplus_{v \in \Si} 
H^1(k_v,G) 
@>>> \bigoplus_{v \in \Si} \hyp^1(k_v,M) @>>> \bigoplus_{v \in \Si} 
H^2(k_v,Y)
\end{CD}
$$
\noindent An inspection of the diagram reveals that for the finiteness of
$D^1(U,M)\{\ell\}$ it suffices to show (using finiteness of $\Si$) the 
finiteness of the torsion groups
$D^1(U,G)\{\ell\}$,  $H^1(k_v,Y)$ and $D^2(U, Y)$, respectively.

\noindent For the first, note that the assumption on $\Sha^1(A_k)$ implies
that $D^1(U,A)\{\ell\}$ is finite by \cite{adt}, \Romannumeral2.5.5, whence the
required finiteness follows from the finiteness of
$H^1(U,T)$ (\cite{adt}, \Romannumeral2.4.6). We have seen the second
finiteness several times when discussing local duality. The finiteness of $D^2(U, Y)$ follows from that of $H^2_c(U,Y)$ which is
is dual to the finite group $H^1(U,T)$  by a result of Deninger (\cite{adt},
\Romannumeral2.4.6). (One can also
prove the finiteness of $D^2(U, Y)$ by the following more direct reasoning: using a 
restriction-corestriction argument, 
one reduces to the case $Y=\Z$. Then $D^2(U,\Z)=D^1(U,\Q/\Z)$ is the dual of
the Galois group of the maximal abelian extension of $k$ unramified
over $U$ and totally split at the places outside $U$, and hence is finite by global class field theory.)

\enddem

\begin{rema} \label{finsha}
{\rm 
The same argument shows that the finiteness of $\Sha^1(A_k)$ implies 
the finiteness of $D^1(U,M)$.
}
\end{rema}

\section{Global results: relation to Galois cohomology}\label{foursect} 

The notation and assumptions are the same as in the previous section.
In particular $U$ is an affine open subset of $\spec \calo_k$ and $\Si
\subset \Omega_k$ consists of the real places and the finite places 
of $\spec \calo_k \setminus U$. Fix an algebraic closure $\kbar$ of $k$ 
(corresponding to a geometric point $\bar \eta$ of $U$).
Let $\Gamma_{\Si}=\pi_1(U, \bar \eta)$ be the 
Galois group of the maximal subfield $k_{\Si}$ 
of $\kbar$ such that the extension $k_{\Si}/k$ is unramified outside 
$\Si$. When $G_k$ is the restriction of a $U$-group scheme $G$ to 
$\spec k $, we shall write 
$H^i(\Ga_{\Si},G_k)$ for $H^i(\Ga_{\Si},G_k(k_{\Si}))$.

\smallskip

We begin with the following analogue of \cite{adt}, \Romannumeral2.2.9:

\begin{prop} \label{firstsha}
Let  $M=[Y \ra G]$ be a 1-motive over $U$ and $M_k$ its restriction 
to $\spec k$. 
Let $\ell$ be a prime number invertible on $U$. Then:

\begin{enumerate}

\item The natural map $\hyp^i(U,M)\{ \ell \} \ra \hyp^i(\Ga_{\Si},M_k)\{ \ell \}$ 
is an isomorphism for $i \geq 1$.

\item The natural 
map $\hyp^1(U,M)\{ \ell \} \ra \hyp^1(k,M_k)\{ \ell \}$ is injective. 

\item The natural map $\hyp^0(U,M) \ra \hyp^0(k,M_k)$ is injective.

\end{enumerate}

\end{prop}

\dem To prove (1) in the case $i>1$, pass to the direct limit over powers of $\ell$ in exact sequence (\ref{nmultiply}). Then Lemma \ref{sorites} (1) implies that it is sufficient to prove the statement for $M$ replaced by $T(M)\{\ell\}$, which follows from \cite{adt}, \Romannumeral2.2.9 by passing to the limit. In the case $i=1$ one unscrews $M$ using the weight filtration and applies \cite{adt}, \Romannumeral2.2.9, \Romannumeral2.4.14 and \Romannumeral2.5.5., 
using also the fact that $A(k)=H^0(U,A)=H^0(\Ga_{\Si},A_k)$ for the abelian 
scheme $A$ with generic fibre $A_k$.
\smallskip

\noindent For (2), it is sufficient by (1) to show that the canonical map 
$\hyp^1(\Ga_{\Si},M_k) \ra \hyp^1(k,M_k)$ is injective. Again using the weight 
filtration on $M$, one reduces this to the injectivity of 
$\hyp^1(\Ga_{\Si}, G_k) \ra \hyp^1(k,G_k)$, 
the injectivity of $H^2(\Ga_{\Si}, Y_k)
\rightarrow H^2(k, Y_k)$ and the surjectivity of $H^1(\Ga_{\Si}, Y_k) 
\ra H^1(k,Y_k)$ (Note that $Y_k(\kbar)=Y_k(k_{\Si})$ because there exists 
a finite \'etale covering $\widetilde U/U$ such that $Y \times_U 
\widetilde U$ is constant).
The two injectivities are consequences of the 
restriction-inflation sequence for $H^1$ in Galois cohomology, 
noting the isomorphisms $$H^2(\Ga_{\Si}, Y_k)\cong H^1(\Ga_{\Si},Y_k(\kbar) 
\otimes \Q/\Z)\quad\hbox{and}\quad H^2(k, Y_k)\cong H^1(k,Y_k(\kbar) 
\otimes \Q/\Z)$$ for the lattice $Y_k(\kbar)=Y_k(k_{\Si})$.
The surjectivity follows from the triviality of the abelian 
group $H^1(\gal(\kbar/k_{\Si}),
Y_k(\kbar))=\Hom_{\rm cont}(\gal(\kbar/k_{\Si}),Y_k(\kbar))$, the Galois group
$\gal(\kbar/k_{\Si})$ being profinite and $Y_k(\kbar)$ torsion-free. 

\smallskip

\noindent The statement (3) is obvious for $M=[0 \ra G]$ (a morphism
from $U$ to $G$ is trivial if and only if it is trivial at the generic
point). Using the weight filtration,
it is sufficient to prove that the map $H^0(U,Y) \ra H^0(k,Y_k)$
is surjective and the map $H^1(U,Y) \ra H^1(k,Y_k)$ is injective.
Actually $H^0(U,Y)=H^0(k,Y_k)$ and $H^0(U,Y/nY)=H^0(k,Y_k/nY_k)$ for each
$n >0$ because $Y$ and $Y/nY$ are locally constant in the \'etale topology.
The injectivity of the map
$H^1(U,Y) \ra H^1(k,Y_k)$ now follows from the commutative exact diagram 
$$
\begin{CD}
H^0(U,Y) @>>> H^0(U,Y/nY) @>>> H^1(U,Y)[n] @>>> 0\\
@VVV @VVV @VVV\\
H^0(k,Y_k) @>>> H^0(k,Y_k/nY_k) @>>> H^1(k,Y_k)[n] @>>> 0.
\end{CD}
$$
\enddem

\medskip

For a 1-motive $M_k$ over $k$ and $i \geq 0$, define the 
{\em Tate-Shafarevich groups}
$$\Sha^i(M_k)={\rm Ker} \, [\hyp^i(k,M_k) \ra \prod_{v \in \Omega_k}
\hyp^i({\hat k}_v,M_k)]$$
If  $M_k$ is the restriction of a 1-motive $M$ defined over $U$, we also define
$$\Sha^i_{\Si}(M_k):={\rm Ker} \, [\hyp^i(\Ga_{\Si},M_k) \ra
\prod_{v \in \Si} \hyp^i({\hat k}_v,M_k)]$$

\begin{rema} \label{hensrem}
{\rm By Lemma~\ref{green}, we can replace ${\hat k}_v$ by $k_v$ in the 
definition of $\Sha^i$ for $i \geq 0$.
}
\end{rema}

Proposition~\ref{firstsha}, Remark~\ref{hensrem} and Corollary \ref{shacor} imply immediately:

\begin{cor}
Let $M$ be a 1-motive over $U$ and $\ell$ a prime number invertible on 
$U$. Then the pairing 
(\ref{interm}) can be identified with a pairing
$$\Sha^1_{\Si}(M_k)\{ \ell \} \times \Sha^1_{\Si}(M^{\ast} _k)\{ \ell \}
\ra \Q/\Z,$$
nondegenerate modulo divisible subgroups.
\end{cor}

In the remaining of this section, we prove Theorem \ref{second}. As in (\cite{adt}, \Romannumeral2.6), the idea is to identify the Tate-Shafarevich groups with the groups $D^i(U, M)$ considered in the previous section for $U$ sufficiently small. But a difficulty is that 
when $Y$ is not trivial, the restriction of $D^i(U,M)$ to $\hyp^i(V,M)$
for $V \subset U$ need not be a subset of $D^i(V,M)$, as the following
example shows. 

\begin{ex}\rm
Let $k$ be a totally imaginary number field with $\pic \calo_k \neq 0$ and 
$M:=\Z[1]$. Then for $U=\spec \calo_k$ we have $D^1(U,M)=H^1(U, \Q/\Z)\neq 0$
by global class field theory. Let $\alpha \neq 0$ in 
$H^1(U, \Q/\Z)$. Then there is a finite place $v$ coming from a closed point
of $U$ such that  the restriction of $\alpha$ to $H^1(k_v,\Q/\Z)$ is nonzero, for otherwise we
would get $\alpha=0$ by Chebotarev's density theorem. Therefore the restriction of $\alpha$ to $H^1(V,M)$ 
does not belong to $D^1(V,M)$ for $V=U-\{v \}$.
\end{ex}

Thus the identification of $\Sha^1(M_k)$ with some $D^1(U,M)$ for suitable $U$
is not straightforward, in contrast to the case of abelian varieties. 
The following proposition takes care of this problem.

Note that by Proposition~\ref{firstsha}, 
any $D^1(U, M)\{ \ell \}$ is naturally a subgroup of $H^1(k, M_k)$ for
$\ell$ invertible on $U$. 

\begin{prop} \label{restriction}
Let $M=[Y \ra G]$ be a 1-motive over $U$ and $\ell$ a prime invertible on $U$. There exists an open subset 
$U_0 \subset U$ such that for any open subset $U_1 \subset U_0$, 
the group $D^1(U_1, M)\{\ell\}$ as a subgroup of $H^1(k, M_k)$ is contained in  $\Sha^1(M_k)\{\ell\}$. 
\end{prop}

For the proof of this proposition we need two preliminary lemmas.

\begin{lem}\label{inj} Let $v$ be a finite place of $U$ and $\calo_v$ the 
henselisation of the local ring of $U$ at $v$. Then the natural map $H^2(\calo_v,Y)\to H^2(k_v,Y)$ is injective when restricted to $\ell$-primary torsion.
\end{lem}

Note that the groups occurring in the lemma are torsion groups (for the first one, this follows from Lemma \ref{sorites} (1) applied to $M=Y[1]$ and passing to the limit).

\begin{dem} By the localisation sequence for the pair $\spec k_v\subset\spec\calo_v$ it is enough to prove triviality of the group $H^2_v(\calo_v, Y)[\ell^n]$ for any $n$. This group is a quotient of $H^1_v(\calo_v, Y/\ell^nY)$ which, according to \cite{adt}, \Romannumeral2.1.10 (a), is a finite group dual to $H^2(\calo_v, \underline{Hom}(Y/\ell^nY, \G))$. But for any finite sheaf $F$ over $\calo_v$ of order prime to the characteristic of the residue field $\k(v)$ of $v$, we have
$H^2(\calo_v,F)=0$ 
because $H^2(\calo_v,F)=H^2(\k(v), \widetilde F)$ by \cite{milne}, 
\Romannumeral3.3.11 a) (where $\widetilde F$ is the restriction of $F$
to $\spec \k(v)$), and $\k(v)$ is of cohomological dimension 1. 
\end{dem}

\begin{lem}\label{sha2}
There exists an open subset 
$U_0 \subset U$ such that for any open subset $U_1 \subset U_0$, 
the group $\Sha^2(Y_k)\{ \ell \}$ contains $D^2(U_1,Y)\{ \ell \}$.
\end{lem}

Note that according to Proposition \ref{firstsha}, for any open $V\subset U$ the group $D^2(V, Y)\{ \ell \}
=D^1(V,Y[1])\{ \ell \}$ identifies with a subgroup of $H^2(k, Y)$. 

\begin{dem} By definition we have
$\Sha^2(Y_k)\{ \ell \} \supset \bigcap_{V \subset U} D^2(V,Y)\{ \ell \}$.
The group $D^2(V,Y)$ is finite for each open subset
$V \subset U$ (cf. proof of Proposition~\ref{partsha}).
Therefore there exist finitely many open subsets $V_1,...,V_r$
of $U$ such that $\Sha^2(Y_k)\{ \ell \} \supset \bigcap_{i=1} ^r
D^2(V_i,Y) \{ \ell \}$.
Let $U_0$ be the intersection
$\bigcap_{i=1}^r V_i$. As the groups $H^2_c(\_\_, Y)$ are covariantly functorial for open immersions $U_1\subset U_2$ (see the previous section), for $i=1,...,r$ there are natural maps $H^2_c(U_0, Y)\to H^2_c(V_i, Y)\to H^2(k, Y)$ which factor through $D^2(U_0, Y)$, so that we get inclusions $D^2(U_0, Y)
\hookrightarrow D^2(V_i, Y)$ and finally 
$D^2(U_0,Y) \{ \ell \} \hookrightarrow\Sha^2(k, Y)\{ \ell \}$.
The above of course holds for any $U_1 \subset U_0$ instead of $U_0$.
\end{dem}

\noindent{\bf Proof of Proposition \ref{restriction}:} Choose $U_0$ as in the previous lemma and take $U_1\subset U_0$. It suffices to show that for each closed point $v\in U_1$ the group $D^1(U_1, M)\{\ell\}$ maps to 0 by the restriction map $\hyp^1(U_1, M)\to\hyp^1(k_v, M)$. Now there is a 
commutative diagram with exact rows

$$
\begin{CD}
H^1(U_1,G) @>>> \hyp^1(U_1,M) @>>> H^2(U_1,Y) \\
@VVV @VVV @VVV \\
H^1(\calo_v,G) @>>> \hyp^1(\calo_v,M) @>{\alpha}>> H^2(\calo_v,Y) \\
@VVV @VVV @VV{\beta}V \\
H^1(k_v,G) @>>> \hyp^1(k_v,M) @>{\gamma}>> H^2(k_v,Y)
\end{CD}
$$
Here the group 
$H^1(\calo_v,G)$ is zero because $G$ is smooth and connected 
over $U_1$, hence $H^1(\calo_v,G)=H^1(\F_v,\widetilde G)$ (cf. \cite{milne},
\Romannumeral3.3.11 a))
is trivial by Lang's Theorem. So the map $\alpha$ in the diagram is injective, but when restricted to $\ell$-primary torsion, the map $\beta$ is injective as well, by Lemma \ref{inj}. Therefore the image of the map $D^1(U_1, M)\{\ell\}\to\hyp^1(k_v, M)$ injects into $H^2(k_v, Y)$ by $\gamma$, so we reduce to the triviality of the composite map $D^1(U_1, M)\{\ell\}\to H^2(k_v, Y)$. This in turn follows from Lemma \ref{sha2} as the map factors through $D^2(U_1, Y)\{\ell\}$ by the right half of the diagram.\enddem

We are now able to prove one half of Theorem \ref{second}:

\begin{theo} \label{mainsha}
Let $M_k$ be a 1-motive over $k$. Then there exists a canonical pairing
$$\Sha^1(M_k) \times \Sha^1(M_k ^{\ast}) \ra \Q/\Z$$
whose kernels are the maximal divisible subgroups of each group.
\end{theo}

\dem We construct the pairing separately for each prime $\ell$. Let $U$ be an open subset of $\spec \calo_k$ such that $M_k$ 
is the restriction of a 1-motive $M$ defined over $U$ and $\ell$ is invertible on $U$. Take a subset $U_0$ as in Proposition~\ref{restriction}. We contend that the inclusion $D^1(U_0, M)\{\ell\}\subset \Sha^1(M_k)\{\ell\}$ furnished by the proposition is in fact an equality, from which the theorem will follow by Corollary \ref{shacor}. Indeed, any element of $\Sha^1(M_k)\{\ell\}$ is contained in some $D^1(V, M)\{\ell\}$, where we may assume $V\subset U_0$. But $D^1(V, M)\{\ell\}\subset D^1(U_0, M)\{\ell\}\subset\Sha^1(M_k)$, by the same argument as in the end of the proof of Lemma \ref{sha2}.

\enddem

\begin{cor} \label{finalcor}
Let $M_k$ be a 1-motive over $k$. 
If $\Sha^1(A_k)$ and $\Sha^1(A_k^{\ast})$ are finite, 
then there is a perfect pairing of finite groups
$$\Sha^1(M_k) \times \Sha^1(M_k ^{\ast}) \ra \Q/\Z$$
\end{cor}

\dem Apply Theorem~\ref{mainsha}, Proposition~\ref{partsha}, and 
Remark~\ref{finsha}.
\enddem

\begin{rema}
{\rm In the case $M_k=[0 \ra T_k]$, we recover the duality between 
$\Sha^1(T_k)$ and $\Sha^2(Y_k^{\ast})$, where $Y_k^{\ast}$ is the module 
of characters of the torus $T_k$. See \cite{tatenak} and 
\cite{neukirch}, \Romannumeral8.6.8.
}
\end{rema}

Now to the result for $\Sha^0$ and $\Sha^2$. First we show:

\begin{lem} \label{zerofini}
Let $M_k$ be a 1-motive over $k$. Then $\Sha^0(M_k)$ is finite.
\end{lem}

\dem We prove that the map $\hyp^0(k, M_k)\to \hyp^0(\hat k_v, M_k)$ has finite kernel
for any finite place $v$ of $k$. With the same notation as above, 
there is a commutative diagram with exact rows:
$$
\begin{CD}
H^0(k,Y_k) @>>> H^0(k,G_k) @>>> \hyp^0(k,M_k) @>>> H^1(k,Y_k) \\
@VVV @VVV @VVV @VVV \\
H^0(\hat k_v,Y_k) @>>> H^0(\hat k_v,G_k) @>>> \hyp^0(\hat k_v,M_k) 
@>>> H^1(\hat k_v,Y_k)
\end{CD}
$$
The group $H^1(k,Y_k)$ is finite because $Y_k(\kbar)$ is a lattice. This being said, a diagram chase using the injectivity of the second vertical map reveals that it suffices to establish the finiteness of the cokernel of the first vertical map. Now this cokernel is certainly of finite type (as $H^0(\hat k_v,Y_k)$ is), and it is annihilated by the degree of  a finite  field extension trivialising the Galois action on $Y_k(\kbar)$, by a
restriction-corestriction argument. 
\enddem

\begin{prop} \label{shazero}
Let $M_k$ be a 1-motive over $k$. Then there is 
a canonical pairing
$$\Sha^0(M_k) \times \Sha^2(M_k ^{\ast}) \ra \Q/\Z$$
whose left kernel is trivial and right kernel is the maximal divisible 
subgroup of $\Sha^2(M_k ^{\ast})$. 
\end{prop}

\dem By Lemma~\ref{zerofini}, it is sufficient to 
construct (for each prime number $\ell$) a pairing 
$\Sha^0(M_k)\{ \ell \} \times \Sha^2(M_k ^{\ast})\{ \ell \} \ra \Q/\Z$
with the required properties.
Let $U$ be an open subset of $\spec \calo_k$ such that $\ell$
is invertible on $U$ and $M$ extends to a 1-motive
$M$ over $U$. Since $D^0(V,M) \{ \ell \}$ is finite for each open subset 
$V$ of $U$ by lemma \ref{sorites} (3), the same argument as in Lemma~\ref{sha2} shows that 
there exists $U_0 \subset U$ such that $\Sha^0(M_k)\{\ell\}=D^0(V,M)\{ \ell \}$ 
for any open subset $V$ of $U_0$ (recall that $D^0(V,M)$ injects into
$\hyp^0(k,M_k)$ by Lemma~\ref{firstsha}, (3)). 

Now notice that $\hyp^2(\calo_v,M)=0$ for each finite place $v$ of $U$. 
Indeed, we have $\hyp^2(\calo_v,M)=\hyp^2(\F_v,\widetilde M)$ by \cite{milne},
\Romannumeral3.3.11, and
we have seen in the proof of Theorem~\ref{unram} that 
$H^2(\F_v,{\widetilde G})=0$ 
and $H^3(\F_v,\widetilde Y)$ is zero because $\F_v$ is of strict
cohomological dimension 2, whence the claim by d\'evissage. 
It follows that for each open subset $V$ of $U$, 
the image of the restriction map $D^2(U,M) \ra \hyp^2(V,M)$ lies in $D^2(V,M)$.
Using the covariant functoriality of compact support cohomology for open immersions, 
we also have ${\cal D}^2(V,M) \subset {\cal D}^2(U,M)$, and finally 
${\cal D}^2(V,M)={\cal D}^2(U,M)$,
where ${\cal D}^2(\_\_, M)$ stands for the image of $D^2(\_\_, M)$ in 
$\hyp^2(k,M_k)$. In particular, $\Sha^2(M_k)={\cal D}^2(V,M)$ for 
each open subset $V$ of $U_0$. 
On the other hand the map $D^2(V,M)\{\ell\} \ra 
\hyp^2(\Gamma_{\Sigma_V},M_k)$ is injective, where $\Sigma_V$ consists of the 
places of $k$ which are not in $V$ (Proposition~\ref{firstsha}, (1)).
Therefore the natural map
$$\underset{^{\longrightarrow}}{{\rm lim}} \,
D^2(V,M)\{ \ell\} \ra \Sha^2(M_k)\{ \ell \}$$
(here the limit is taken over nonempty open subsets $V$ of $U_0$) 
is an isomorphism. The proposition follows from the identification 
of $D^0(V,M)\{ \ell \}$ with $\Sha^0(M_k)\{\ell\}$ and from
the case $i=0$ in Corollary~\ref{shacor}.

\enddem

\begin{remas}\label{shazeroremas}${}$
{\rm 
\begin{enumerate}
\item The special case $M_k=Y_k[1]$ corresponds to the duality between 
the groups 
$\Sha^1(Y_k)$ and $\Sha^2(T_k^{\ast})$, where $T_k^{\ast}$ is the torus 
with character module $Y_k$. Compare \cite{adt}, \Romannumeral1.4.20 (a).
\item The case $M_k=[0\to A_k]$ does not contain much new information. Indeed, assume for instance $k$ totally imaginary. Then combining the injectivity of 
$A_k(k)\to A_k({\hat k}_v)$ as well as the vanishing of $H^2({\hat k}_v, A_k)$ 
for each finite place $v$ (\cite{cogal}, \Romannumeral2.5.3, Proposition 16) with the above proposition we get that $H^2(k, A_k)$ is a divisible group; 
but this we know anyway since $k$ is of cohomological dimension 2. In fact 
in this case $H^2(k, A_k)=0$ by \cite{adt}, \Romannumeral1.6.24.
\item Presumably, $\Sha^2(M)$ is finite, as it is finite for $M=T$ and even trivial for $M=A$ and $M=Y[1]$.
\end{enumerate}
}
\end{remas}

\section{The Poitou-Tate exact sequence} \label{fivesect}

We keep notation from the previous sections. In 
particular, $k$ is a number field, $M_k$ is a 1-motive over $k$, and $U$ is 
an open  subset of $\spec \calo_k$ such that $M_k$ extends to 
a 1-motive $M$ over $U$.
For each finite place $v$ of $k$, we denote by $\widehat \calo_v$ the ring of 
integers of the completion ${\hat k}_v$. The groups $\hyp^i(k,M_k)$ 
($-1 \leq i \leq 2$) and $\hyp^i(\hat k_v, M_k)$ for $i \geq 1$ 
($v \in \Omega_k$) are equipped with the discrete topology. 
For $i \geq 0$, we define $\P^i(M_k)$ as the restricted product
over $v \in \Omega_k$ of the $\hyp^i({\hat k}_v,M_k)$, with respect to 
the images of $\hyp^i(\widehat \calo_v,M)$ for finite places $v$ of $U$.
The groups $\P^i(F_k)$ are defined similarly for any \'etale sheaf 
$F_k$ over $\spec k$. 
Clearly all these groups are independent of the choice of $U$. We equip them
with their restricted product topology. 
Note that $\P^2(M_k)$ is the direct sum of the groups $\hyp^2(\hat k_v,M_k)$ for all places $v$:
indeed, for each finite place $v$ of $U$, we have $\hyp^2(\hat \calo_v,M)=\hyp^2(\F_v,\widetilde M)$, and we have already seen that this last group vanishes.
We have natural restriction maps $\beta_i:\,\hyp^i(k, M_k)\to \P^i(M_k)$ for
$i=1,2$; their kernels are precisely the groups $\Sha^i(M_k)$. 

In this section we establish a Poitou-Tate type exact sequence for $M_k$. A technical complication for our considerations to come arises from the fact that one is forced to work with two kinds of completions: the first one is what we have used up till now, i.e. the inverse limit $A^{\wedge}$ of all open subgroups of finite index of a topological abelian group $A$, and the second is the inverse limit $A_{\wedge}$ of the quotients $A/n$ for all $n>0$. These are not the same in general (indeed, the latter is not necessarily profinite). But it is easy to see that $(A_{\wedge})^{\wedge}$ is naturally isomorphic to $A^{\wedge}$; in particular, we have a natural map $A_{\wedge}\to A^{\wedge}$.  

The natural map $\hyp^0(k, M)\to \P^0(k, M)$ therefore gives rise to a commutative diagram between the two different kinds of completions:
\begin{equation}\label{compl}
\begin{CD}
\hyp^0(k, M)_{\wedge} @>\theta_0>> \P^0(k, M)_{\wedge}\cr
@VVV @VVV \cr
\hyp^0(k, M)^{\wedge} @>\beta_0>> \P^0(k, M)^{\wedge}.
\end{CD}
\end{equation}

Denote by $\Sha^0_{\wedge}(M_k)$ the kernel of the above map $\theta_0$. We first prove the following duality result analogous to \cite{adt}, \Romannumeral1.6.13 (b). 

\begin{prop} \label{zerotheo}
Let $M_k=[Y_k\to G_k]$ be a 1-motive over $k$. Assume that $\Sha^1(A_k)$ is finite. Then there is a perfect 
pairing 
$$\Sha^0_{\wedge}(M_k) \times \Sha^2(M_k ^{\ast}) \ra \Q/\Z$$
where the first group is compact and the second is discrete.
\end{prop}

\begin{remas}\rm${}$
\begin{enumerate}
\item The proof below will show that if one only assumes the finiteness of the $\ell$-primary torsion part of $\Sha^1(A_k)$, then one gets a similar duality between the $\ell$-primary torsion part of $\Sha^2(M_k^{\ast})$ and a group $\Sha_{(\ell)}^0(M_k)$ defined similarly to $\Sha^0_{\wedge}(M_k)$ but taking only $\ell$-adic completions.
\item It is natural to ask how the above proposition is related to 
Proposition \ref{shazero}. In the case when $M_k^{\ast}=[0\to T_k]$ is a torus, 
one knows that $\Sha^2(M_k ^{\ast})=\Sha^2(T_k)$ is finite and so by 
comparing the 
two statements we get that $\Sha^0_{\wedge}(M_k)=\Sha^0(M_k)$ and the two 
dualities are the same. 
In the case $M_k=[0\to A_k]$, both $\Sha^0_{\wedge}(A_k)$ and $\Sha^2(A_k)$
are trivial (\cite{adt}, Corollaries \Romannumeral1.6.23 and 
\Romannumeral1.6.24), but even the fact that they are dual groups 
does not seem to be an easy consequence of Proposition \ref{shazero}. The equality $\Sha^0(M)=\Sha^0_{\wedge}(M)$ would follow via Proposition \ref{shazero} from the finiteness of $\Sha^2(M^{\ast})$, but we do not know this latter fact (compare Remark \ref{shazeroremas} (3)).
\end{enumerate}
\end{remas}

For the proof of the proposition we first show a lemma.

\begin{lem}\label{resprod}
For any $n>0$, the natural map
$$
\P^0(M_k)/n\to\prod_v\hyp^0(\hat k_v, M_k)/n
$$
is injective and its image is the restricted product of the groups $\hyp^0(\hat k_v, M_k)/n$ with respect to the subgroups $\hyp^0(\hat \calo_v, M_k)/n$.\end{lem}

\begin{dem}
Take an element $x=(x_v) \in \P^0(M_k)$ and assume $x_v$ is in $n\hyp^0(\hat k_v,M)$ for all $v$. For almost all $v$, it also comes from $\hyp^0(\hat \calo_v,M)$, hence  in fact it comes from $n\hyp^0(\hat\calo_v,M)$, for in the exact commutative 
diagram 
$$
\begin{CD}
\hyp^0(\widehat \calo_v, M) @>{n}>>\hyp^0(\widehat \calo_v, M) 
@>>> H^1(\widehat \calo_v, T_{\Z/n\Z}(M)) \cr
@VVV @VVV @VVV \cr
\hyp^0(\hat k_v, M_k) @>{n}>> \hyp^0(\hat k_v, M_k) @>>> H^1(\hat k_v,T_{\Z/n\Z}(M_k))
\end{CD}
$$
the third vertical map is injective, $H^1_v(\widehat\calo_v,T_{\Z/n\Z}(M))$ being trivial by the same argument as in the proof of Lemma \ref{inj}. But this means $x\in n\P^0(M_k)$. 
The second statement is obvious.
\end{dem}

\noindent{\bf Proof of Proposition \ref{zerotheo}:}
Using exact sequence (\ref{nmultiply}) we get a commutative diagram 
with exact rows for each $n>0$:
$$
\begin{CD}
0 @>>> \hyp^0(k,M_k)/n @>>> H^1(k,T_{\Z/n\Z}(M_k)) @>>> \hyp^1(k,M_k)[n] 
@>>> 0 \\
&& @VVV @VVV @VVV \\
0 @>>> \hyp^0({\hat k}_v,M)/n @>>> H^1({\hat k}_v,T_{\Z/n\Z}(M_k)) 
@>>> \hyp^1({\hat k}_v,M_k)[n] @>>> 0 
\end{CD}
$$
and a similar exact diagram holds with $\widehat \calo_v$ instead of 
${\hat k}_v$. Taking restricted products and using the above lemma, we get a commutative exact diagram:
$$
\begin{CD}
0 @>>> \hyp^0(k,M_k)/n @>>> H^1(k,T_{\Z/n\Z}(M_k)) @>>> \hyp^1(k,M_k)[n]
@>>> 0 \\
&& @VVV @VVV @VVV \\
0 @>>> \P^0(M_k)/n @>>> \P^1(T_{\Z/n\Z}(M_k))
@>>> \P^1(M_k)[n] @>>> 0
\end{CD}
$$
If we pass to the inverse limit over all $n$, the
two lines remain left exact, whence a commutative exact diagram 
\begin{equation}\label{beta1}
\begin{CD}
0 @>>> \hyp^0(k,M_k)_{\wedge} @>>> H^1(k,T(M_k)) @>>> I_1
@>>> 0 \\
&& @VV{\theta_0}V @VV{\theta}V @VV{\beta_1}V \\
0 @>>> \P^0(M_k)_{\wedge} @>>> \P^1(T(M_k))
@>>> I_2 @>>> 0
\end{CD}
\end{equation}
where we define $\P^1(T(M_k))$ (resp. $H^1(k,T(M_k))$) as the inverse limit 
of $\P^1(T_{\Z/n\Z}(M_k))$ (resp. $H^1(k,T_{\Z/n\Z}(M_k))$).
Here $I_1$, $I_2$ are respectively subgroups of the full Tate modules
$T(H^1(k,M_k))$ and $T(\P^1(M_k))$.
In particular, ${\rm Ker} \, \beta_1$ is a subgroup of 
$T(\Sha^1(M_k))$ because the inverse limit functor is left exact.
But $T(\Sha^1(M_k))$ is zero thanks to the finiteness assumption on 
$\Sha^1(A_k)$ (which implies the finiteness 
of $\Sha^1(M_k)$ as in Corollary~\ref{finalcor}).
Therefore we obtain that $\Sha^0_{\wedge}(M_k)={\rm Ker} \, \theta_0$
is isomorphic to ${\rm Ker} \, \theta$. By Poitou-Tate duality for finite modules (\cite{neukirch}, \Romannumeral8.6.8), there is a perfect pairing between the latter group and
$\Sha^2(T(M_k ^{\ast})_{\tors})$. 
We conclude by observing that 
$\Sha^2(T(M_k ^{\ast})_{\tors})$ is also $\Sha^2(M_k ^{\ast})_{\tors}$
by the same argument as in \cite{adt}, \Romannumeral1.6.8 (use the 
analogue of the exact sequence (\ref{nmultiply}) for $M_k$ over $\spec k$ with $i=2$ and pass to the limit). 
\enddem\medskip

We now return to diagram (\ref{compl}) and prove:

\begin{prop}\label{ker} Keep the finiteness assumption on $\Sha^1(A_k)$.
Then, with notations as in diagram (\ref{compl}), the natural map $\ker\theta_0\nobreak\to\nobreak\ker\beta_0$ is an isomorphism.
\end{prop}

By virtue of the proposition, we may employ the notation $\Sha^0_{\wedge}(M_k)$ for $\ker\beta_0$ as well and use the resulting duality. The finiteness assumption on $\Sha^1(A_k)$ is presumably superfluous here but we did not succeed in removing it (and use it elsewhere anyway).

For the proof we need a lemma about abelian groups.

\begin{lem}
Let $A$ be a discrete abelian group of finite exponent $n$. Then the intersection of the finite index subroups of $A$ is trivial.
\end{lem}

\begin{dem}
Consider the profinite group $A^D=\Hom(A, \Z/n)$. As $A^{DD}=A$, the statement is equivalent to saying that any character of $A^D$ vanishing on all finite subgroups of $A^D$ is trivial. This holds because all finitely generated subgroups of $A^D$ are finite.
\end{dem}

{\noindent \bf Proof of Proposition \ref{ker}:} We begin by showing that the vertical maps in diagram (\ref{compl}) above are injective, whence the injectivity of the map $\ker\theta_0\to\ker\beta_0$. For the left one, note that the topology on $\hyp^0(k, M)$ being discrete, injectivity means that any element of $\hyp^0(k, M)$ which is nontrivial modulo $n$ for some $n$ gives a nonzero element in some finite quotient -- this holds by the lemma above.

For injectivity of the map $\P^0(k, M)_{\wedge}\to\P^0(k, M)^{\wedge}$, take an element $x=(x_v) \in \P^0(M_k)$ not lying in $n\P^0(M_k)$. We have to find an open subgroup of finite index avoiding $x$. By Lemma \ref{resprod}, there is a local component $x_v$ not lying in $n\hyp^0(\hat k_v,M)$. 
We thus get that for some $v$, our $x$ is not contained in the inverse image in $\P^0(M_k)$ of the subgroup $n\hyp^0(\hat k_v,M)\subset\hyp^0(\hat k_v,M)$, which is open of finite index, by definition of the topology on $\hyp^0(\hat k_v,M)$. 

For the surjectivity of the map $\ker\theta_0\to\ker\beta_0$, remark first that $\ker\theta_0$ is a profinite group, being dual to the torsion group $\Sha^2(M_k)$ by the previous proposition. 
Therefore $(\ker\theta_0)^{\wedge}=\ker\theta_0$, so by completing the exact sequence 
\begin{equation}\label{quot}
0\to\ker\theta_0\to\hyp^0(k, M_k)_{\wedge}\to\im\theta_0\to 0
\end{equation} 
we get an exact sequence 
$\ker\theta_0\to\hyp^0(k, M_k)^{\wedge}\to(\im\theta_0)^{\wedge}\to 0$.  
To conclude, it is thus enough to show the injectivity of the natural map $(\im\theta_0)^{\wedge}\to\P^0(M_k)^{\wedge}$. By the above considerations, $\P^0(M_k)_{\wedge}$ injects into its completion $\P^0(M_k)^{\wedge}$, so the completion of the subgroup $\im\theta_0$ is simply its closure in $\P^0(M_k)^{\wedge}$, whence the claim. But there is a subtle point here: in this argument, $\im\theta_0$ is equipped with the subspace topology inherited from $\P^0(M_k)_{\wedge}$, whereas in exact sequence (\ref{quot}) it carried the quotient topology from $\hyp^0(M_k)_{\wedge}$; we have to check that the two topologies are the same. It is enough to show that $\im\theta_0$ is closed in $\P^0(M_k)_{\wedge}$, for then it is locally compact in the subspace topology and \cite{hewitt}, 5.29 applies.   
To do so, observe that the image of the map $\theta$ in diagram (\ref{beta1}) above equals the kernel of the map $\P^1(T(M_k))\to H^1(k, T(M_k^*)_{\rm tors})^D$ coming from local duality, as a consequence of the Poitou-Tate sequence for finite modules (\cite{neukirch}, \Romannumeral8.6.13). Now a diagram chase in (\ref{beta1}) using the injectivity of $\beta_1$ reveals that $\im\theta_0$ is the kernel of the (continuous) composite map $\P^0(M_k)_{\wedge}\to\P^1(T(M_k))\to H^1(k, T(M_k^*)_{\rm tors})^D$, hence it is closed indeed. 
\enddem

We can now state the main result of this section.

\begin{theo}[Poitou-Tate exact sequence] \label{poitsuite}
Let $M_k$ be a 1-motive over $k$. Assume that $\Sha^1(A_k)$ and 
$\Sha^1(A^{\ast}_k)$ are finite, where $A_k$ is the abelian variety 
corresponding to $M_k$.
Then there is a twelve term exact sequence of topological groups

\begin{equation} \label{poitou}
\begin{CD}
0 @>>> \hyp^{-1}(k,M_k)^{\wedge} @>{\gamma_2}^D>> \prod_{v \in \Omega_k} 
\hyp^2(\hat k_v,M^{\ast}_k)^D @>{\beta_2}^D>> \hyp^2(k,M_k^{\ast})^D \\
&&&&&& @VVV \\
&& \hyp^1(k,M^{\ast}_k)^D @<{\gamma_0}<< \P^0(M_k)^{\wedge} @<{\beta_0}<< 
\hyp^0(k,M_k)^{\wedge} \\
&& @VVV \\
&& \hyp^1(k,M_k) @>{\beta_1}>> \P^1(M_k)_{\rm tors} @>{\gamma_1}>> 
(\hyp^0(k,M^{\ast}_k)^D)_{\rm tors} \\
&&&&&& @VVV \\
0 @<<< \hyp^{-1}(k,M^{\ast}_k)^D @<{\gamma_2}<< \bigoplus_{v \in \Omega_k} 
\hyp^2(\hat k_v,M_k) @<{\beta_2}<< \hyp^2(k,M_k) 
\end{CD}
\end{equation}
where the maps $\beta_i$ are the restriction maps defined at the beginning of this section, the maps $\gamma_i$ are induced by the local duality theorem of Section 2, and the unnamed maps come from the global duality results of Proposition \ref{zerotheo} (completed by Proposition \ref{ker}) and Corollary \ref{finalcor}.
\end{theo}

\begin{remas}\rm ${}$
\begin{enumerate}
\item In the above sequence the group $\P^1(M_k)_{\rm tors}$ is equipped with the discrete topology, and {\em not} the subspace topology from $\P^1(M_k)$. 
\item
The sequence (\ref{poitou}) is completely symmetric in the sense that if we 
replace $M_k$ by $M^{\ast}_k$ and dualise, we obtain exactly the same sequence. 
\item
If we only assume the finiteness of $\Sha^1(A_k)\{ \ell \}$ and 
$\Sha^1(A_k ^{\ast})\{ \ell \}$ for some prime number $\ell$, then 
the analogue of the exact sequence (\ref{poitou}) still holds, with profinite completions replaced by $\ell$-adic 
completions  and the torsion groups involved by their $\ell$-primary part. 
\item Some special cases of the theorem are worth noting. 
For $M_k=[0 \ra A_k]$ we get the ten-term exact sequence of 
\cite{adt}, \Romannumeral1.6.14. 
For $M_k=[0 \ra T_k]$ one can show (see proof of Proposition~\ref{derniere} 
below)
that $\beta_2 ^D$ is 
an isomorphism, hence we obtain a classical (?) nine-term exact sequence; 
the case $M_k=Y_k[1]$ is symmetric (compare \cite{adt}, \Romannumeral1.4.20).
\end{enumerate}
\end{remas}

The proof will use the following lemma.

\begin{lem} \label{torsfin}
Let $M_k$ be a 1-motive over $k$. Then $\hyp^0(k,M_k)_{\tors}$ is finite.
\end{lem}

\dem Let $M'_k=M_k/W_{-2}(M_k)=[Y_k \ra A_k]$. From the exact sequence
$$H^0(k,Y_k) \ra H^0(k,A_k) \ra \hyp^0(k,M'_k) \ra H^1(k,Y_k)$$
we deduce that $\hyp^0(k,M'_k)$ is of finite type because $H^0(k,A_k)$
is of finite type (Mordell-Weil theorem) and $H^1(k,Y_k)$ is finite
(because $Y_k(\kbar)$ is a lattice). There is also an exact sequence
$$\hyp^{-1}(k,M'_k) \ra H^0(k,T_k) \ra \hyp^0(k,M_k) \ra 
\hyp^0(k,M'_k)$$
where $T_k$ is the torus corresponding to $M_k$.
It is therefore sufficient to show that the torsion subgroup of the group $B:=T_k(k)/\im 
\hyp^{-1}(k,M'_k)$ is finite. Choose 
generators $a_1,...,a_r$ of the lattice $\hyp^{-1}(k,M'_k)$, and let 
$b_1,...,b_r$ be their images in $T_k(k)$. We can find an open subset $U$ 
of $\spec \calo_k$ such that $T_k$ extends to a torus $T$ over $U$, and 
$b_i \in H^0(U,T)$ for any $i \in \{ 1,...,r \}$. In this way, 
each $x \in B_{\tors}$ is the image in $B$ of an element $y \in T_k(k)$ for which $y^n \in H^0(U,T)$ for some $n >0$. Let $V/U$ be an
\'etale covering such that $T$ splits over $V$, i.e. it becomes isomorphic to some power $\G^N$. If $L$ denotes the fraction field of $V$, we have $T_k(L)/H^0(V, T)\cong (L^{\times}/H^0(V, \G))^N$ which is naturally a subgroup of the free abelian group $\hbox{\rm Div}(V)^N$; in particular, it has no torsion. Therefore, since $y^n \in H^0(V,T)$ via the inclusion $H^0(U, T)\to H^0(V, T)$ we get that $y \in H^0(V,T)=H^0(V, \G)^N$. Let $H$ be the subgroup of $H^0(V, T)$ generated by the $y$'s; since $H^0(V,\G)$ is of finite type by 
Dirichlet's Unit Theorem, so is $H$. We thus get a surjection from the finitely generated group $H$ to the torsion group $B_{\tors}$, whence the claim.\enddem

\noindent{\bf Proof of Theorem \ref{poitsuite}:} The first line is dual to the
last one, so for its exactness it is enough to show the exactness of the latter. We proceed as in \cite{adt}, 
\Romannumeral1.6.13 (b). For each $n >0$, we have an exact commutative diagram 
(using the exact sequence (\ref{nmultiply}) and the Poitou-Tate sequence 
for the finite module $T_{\Z/n\Z}(M_k)$):

$$
\begin{CD}
&&&& 0 \\
&&&& @VVV \\
&&&& \hyp^0(k,M_k^{\ast})[n] ^D \\
&&&& @VVV \\
H^2(k, T_{\Z/n\Z}(M_k)) @>>> \bigoplus_{v \in \Omega_k} 
H^2(\hat k_v,T_{\Z/n\Z}(M_k)) @>>> H^0(k,T_{\Z/n\Z}(M_k ^{\ast}))^D @>>> 0 \\
@VVV @VVV @VVV \\
\hyp^2(k,M_k)[n] @>>> \bigoplus_{v \in \Omega_k} \hyp^2(\hat k_v,M_k)[n]
@>>> (\hyp^{-1}(k,M_k^{\ast})/n)^D  \\
&&&& @VVV \\
&&&& 0
\end{CD}
$$

By Lemma~\ref{torsfin}, the full Tate module $T
(\hyp^0(k,M_k^{\ast}))$ is trivial.
Therefore taking the inductive limit over all $n$, 
we obtain the exact commutative diagram, where the right vertical map is 
an isomorphism:
$$
\begin{CD}
H^2(k, T(M_k)) @>>> \bigoplus_{v \in \Omega_k} H^2(\hat k_v,T(M_k))
@>>> H^0(k,T(M_k ^{\ast}))^D @>>> 0 \\
@VVV @VVV @VVV \\
\hyp^2(k,M_k) @>>> \bigoplus_{v \in \Omega_k} \hyp^2(\hat k_v,M_k) 
@>>> {{(\hyp^{-1}(k,M_k^{\ast})}^{\wedge})}^D @>>> 0
\end{CD}
$$

We remark that the first two vertical maps are also isomorphisms by the same
argument as at the end of the 
proof of Proposition~\ref{zerotheo}. 
Since the first row is exact by the Poitou-Tate sequence
for finite modules (\cite{neukirch}, \Romannumeral8.6.13), so is the second one.
Therefore the exactness of the  last line of (\ref{poitou}) (and hence that of
the first) follows noting that
the dual of the profinite completion of 
the lattice $\hyp^{-1}(k,M_k^{\ast})$ is the same as the dual of 
$\hyp^{-1}(k,M_k^{\ast})$ itself. 

To prove exactness of the second line, note first that it is none but the profinite completion of the sequence
\begin{equation}\label{sequ}
H^0(k, M_k)_{\wedge}\stackrel{\theta_0}{\to}\P^0(M_k)_{\wedge}\stackrel{\gamma'_0}\to H^1(k, M_k^*)^D, 
\end{equation}
where the map $\gamma'_0$ is induced by local duality, taking Lemma \ref{resprod} into account. First we show that this latter sequence is a complex.
The map $\gamma'_0$ has the following concrete description at a finite level:
for $\alpha=(\alpha_v)\in \P^0(M_k)/n$ 
and $\beta \in \hyp^1(k,M^{\ast}_k)[n]$, denote by $\beta_v$ the image of 
$\beta$ in $\hyp^1(\hat k_v,M^{\ast}_k)[n]$. 
Then $[\gamma'_0(\alpha)](\beta)$ is the sum over all $v$  of the elements
$j_v(\alpha_v\cup\beta_v)$, where $(\alpha_v\cup\beta_v)\in H^2(\hat k_v, \mu_n)$
via the pairing $M\times M^{\ast}\to\G[1]$, and $j_v$ is the local invariant.
(The sum is finite by virtue of the property $H^2(\widehat \calo_v,\mu_n)=0$
because the 
elements $\alpha_v$ and $\beta_v$ are unramified for almost all $v$.) 
Then
$\gamma_0 '\circ \theta_0=0$ follows (after passing to the limit) from the reciprocity law of global class field theory, according to which the sequence
\begin{equation}\label{brauer} 
0\ra H^2(k,\mu_n) \ra 
\bigoplus_{v \in \Omega_k} H^2(\hat k_v,\mu_n) \stackrel{\sum j_v}{\ra} \Z/n\to 0
\end{equation} 
is a complex.

For the exactness of the sequence (\ref{sequ})
recall that, as remarked at the end of the proof of Proposition \ref{ker}, diagram (\ref{beta1}) and the Poitou-Tate sequence for finite modules imply that $\im\theta_0$ is the kernel of the composed map 
$$
\P^0(M_k)_{\wedge}\to\P^1((T(M_k))\to H^1(k, T(M_k^*)_{\rm tors})^D. 
$$
The claimed exactness then follows from the commutative diagram
$$
\begin{CD}
\P^0(M_k)_{\wedge} @>>> \P^1(T(M_k)) \\
 @VV{\gamma'_0}V @VVV \\
\hyp^1(k,M_k^{\ast})^D @>>> \hyp^1(k,T(M_k^{\ast})_{\tors})^D
\end{CD}
$$
whose commutativity arises from the compatibility of the duality pairings for
1-motives and their ``$n$-adic'' realisations via the Kummer map.

Now we show that the profinite completion of sequence (\ref{sequ}), i.e. the second row of diagram (\ref{poitou}) remains exact. This follows by an argument similar to the one at the end of the proof of Proposition \ref{ker}, once having checked that $\im\gamma'_0$ is closed in $\hyp^1(k, M_k^{\ast})^D$ and $(\im \gamma'_0)^{\wedge}=\im \gamma'_0$. To see this, note that by applying the snake lemma to diagram (\ref{beta1}) we get that ${\rm Coker} \, \theta_0$
(with the quotient topology) injects as a closed subgroup into ${\rm Coker} \, \theta$. But one sees using the Poitou-Tate sequence for finite modules that the latter group is profinite, hence so is $\im \gamma'_0$; in particular, it is compact and hence closed in $\hyp^1(k, M_k^{\ast})^D$.

Next, remark that by definition of the restricted product topology 
and Theorem~\ref{unram}, the dual of the group $\P^0(M_k)$
(equipped with the restricted product topology) is $\P^1(M_k^{\ast})$.
Thus the dual of $\P^0(M_k)^{\wedge}$ is $\P^1(M_k^{\ast})_{\rm tors}$.
Therefore we obtain the third line by dualising the second 
one (and exchanging the roles of $M$ and $M^{\ast}$), which consists 
of profinite groups.

Finally, the 
exactness of the sequence (\ref{poitou}) at the ``corners'' follows
immediately, in the first two rows, from the dualities $\Sha^0_{\wedge}(M_k)\cong\Sha^2(M_k ^{\ast})^D$ (Proposition~\ref{zerotheo} combined with Proposition \ref{ker}) and $\Sha^1(M_k)\cong\Sha^1(M_k^{\ast})^D$ (Corollary~\ref{finalcor}), respectively; the remaining corner is dual to the first one.
\enddem

We conclude with the following complement.

\begin{prop} \label{derniere}
Let $M_k$ be a 1-motive over $k$. Then the natural map 
$\hyp^i(k,M_k) \ra \bigoplus_{v \in \Omega_{\infty}} \hyp^i(k_v,M_k)$ is an 
isomorphism for $i \geq 3$.
\end{prop}

\dem When $M=[0 \ra G]$, this follows immediately by devissage from \cite{adt}, 
\Romannumeral1.4.21, and \Romannumeral1.6.13 (c).
To deal with the general case, it is sufficient to show that the 
map $f_i : 
H^i(k,Y_k) \ra \bigoplus_{v \in \Omega_{\RR}} H^i(k_v,Y_k)$ is an 
isomorphism for $i \geq 3$. 
Using the exact sequence 
$$0 \ra H^i(k,Y_k)/n \ra H^i(k,Y_k/nY_k) \ra H^{i+1}(k,Y_k)[n] \ra 0$$
for each $n >0$, we reduce to the case $i=3$. The last line of the Poitou-Tate exact 
sequence for $M_k=Y_k[1]$ yields the surjectivity of $f_3$ because 
$\hat k_v$ is of strict cohomological dimension 2 for $v$ finite. 
On the other hand, we remark that $H^3(k,\Z)=0$ (\cite{adt}, 
\Romannumeral1.4.17), hence $H^3(k,Y_k)=H^3(k,Y_k)[n]$ for some $n >0$ by a 
restriction-corestriction argument. In particular the divisible 
subgroup of $H^3(k,Y_k)$  is zero.
The injectivity of $f_3$ now 
follows from Proposition~\ref{shazero} applied to $M_k=T_k^{\ast}$, where 
$T_k^{\ast}$ is the torus with module of characters $Y_k$.

\enddem

\begin{rema}
{\rm 
The global results in this paper still hold (with the same proofs) 
if we replace
the number field $k$ by a function field in one variable over a 
finite field, provided that one ignores the $p$-primary torsion part of the groups under consideration, where $p={\rm char} \, k$). We leave the verification of this to the readers.
}
\end{rema}

\section{Comparison with the Cassels-Tate pairing}

In this section, we give a definition of the pairing of Theorem \ref{second} purely in terms of Galois cohomology and show that in the case $M=[0\to A]$ it reduces to the classical Cassels-Tate pairing for abelian varieties.

The idea is to use the {\em diminished cup-product} construction discovered by Poonen and Stoll (see \cite{poonen}, pp. 1117--1118). One could present it in a general categorical setting but for the ease of exposition we stick to the special situation we have. So assume given three exact sequences
$$
0\to M_1\to M_2\to M_3\to 0,
$$
$$
0\to N_1\to N_2\to N_3\to 0,
$$
$$
0\to P_1\to P_2\to P_3\to 0
$$
where the $M_i$, $N_i$, $P_i$ are complexes of abelian \'etale sheaves over some scheme $S$. Assume further given a pairing $M_2\otimes^{\bf L} N_2\to P_2$ in the derived category that maps $M_1\otimes^{\bf L} N_1$ to $P_1$. Then one has a pairing
$$
{\rm Ker}\, [\hyp^i(S, M_1)\!\to\!\hyp^i(S, M_2)]\times
{\rm Ker}\,[\hyp^j(S, N_1)\!\to\!\hyp^j(S, N_2)]\!\to\! \hyp^{i+j-1}(S, P_3)
$$
defined as follows. Any element $a$ of $\ker [\hyp^i(S, M_1)\to\hyp^i(S,
M_2)]$ comes from an element $b$ of $\hyp^{i-1}(S, M_3)$. Since the above
pairing induces a pairing $M_3\otimes^{\bf L}N_1\to P_3$, our $b$ can be
cupped with $a'\in \ker [\hyp^j(S, N_1)\to\hyp^j(S, N_2)]$ to get an element
in $\hyp^{i+j-1}(S,P_3)$. It follows from the fact that $a'$ maps to 0 in
$\hyp^j(S, N_2)$ that this definition does not depend on the choice of $b$.

Now apply this in the following situation. Let ${\bf A}_k$ be the ad\`ele ring of the number field $k$, $S=\spec k$ and $M_k$ a 1-motive over $k$. Consider \'etale sheaves over $\spec k$ as $\gal(k)$-modules and put $M_1=M_k(\kbar)$, $M_2=M_k(\kbar\otimes_k{\bf A}_k)$, $N_1=M_k^{\ast}(\kbar)$, $N_2=M_k^{\ast}(\kbar\otimes_k{\bf A}_k)$, $P_1=\G(\kbar)[1]$, $P_2=\G(\kbar\otimes_k{\bf A}_k)[1]$ and define $M_3$, $N_3$, $P_3$ to make the above sequences exact. The pairing $M_k\otimes^{\bf L}M_k^{\ast}\to\G[1]$ induces the pairing $M_2\otimes^{\bf L}N_2\to P_2$ required in the above situation. Note that we have   
$$
\Sha^i(M_k)={\rm Ker}\,[\hyp^i(k, M_1)\to\hyp^i(k, M_2)],$$
$$\Sha^j(M^{\ast}_k)={\rm Ker}\,[\hyp^j(k, N_1)\to\hyp^j(k, N_2)].
$$ 
Finally, class field theory tells us that the cokernel of the map $H^2(k,
\G)\to H^2(k, \G(\kbar\otimes_k {\bf A}_k))$ is isomorphic to $\Q/\Z$; indeed,
one has $$H^2(k, \G(\kbar\otimes_k {\bf A}_k))=\bigoplus_v H^2(k_v, \G)$$ by
Shapiro's lemma (and the isomorphisms $H^2(k_v, \G)\cong H^2(\hat
k_v, \G)$, cf. Lemma \ref{green}), so the claim follows from the exactness of the sequence (\ref{brauer}) (\cite{neukirch}, Theorem 8.1.17) in view of the vanishing of $H^3(k, \G)$ (\cite{neukirch}, 8.3.10 (iv)).

Putting all this together, we get that the diminished cup-product construction yields for $j=2-i$ pairings
\begin{equation}\label{finalpair}
\Sha^i(M_k)\times\Sha^{2-i}(M_k)\to\Q/\Z
\end{equation}  
for $i=0,1$.

\begin{prop}
The above pairings coincide with those constructed in Section 4.
\end{prop}

\begin{dem}
Assume $M_k$ extends to a 1-motive $M$ over an open subset $U\subset \spec k$; let $\Sigma$ denote the finite set of places of $k$ which are real or coming from closed points of $\spec\calo_k$ outside $U$. Apply the diminished cup-product construction for $S=U$, $M_1=M$, $M_2=p_*p^*M$, $N_1=M^{\ast}$, $N_2=p_*p^*M^{\ast}$, $P_1=\G[1]$, $P_2=p_*p^*\G[1]$, where $p:\,\oplus_{v\in\Sigma}\spec k_v\to U$ is the canonical morphism (here we use the conventions of Section 3 for $k_v$). Note that, in the notation of Section 3, 
$$
D^i(U,M)=\ker [\hyp^i(U, M_1)\to\hyp^i(U, M_2)]$$
and
$$
D^j(U,M^{\ast})=\ker [\hyp^j(U, N_1)\to\hyp^j(U, N_2)].
$$ 
Moreover, one has $${\rm coker}\,[H^2(U, \G)\to H^2(U, p_*p^*\G)]\cong H^3_c(U, \G)\cong \Q/\Z$$ where the first isomorphism comes from Shapiro's lemma in the \'etale setting and the second from class field theory, both combined with the fundamental long exact sequence of compact support cohomology (the second isomorphism was already used to construct the pairings of Chapter 3). For $i=0,1$, $j=2-i$ one thus gets pairings
$$
D^i(U,M)\times D^{2-i}(U, M)\to \Q/\Z
$$
which are manifestly the same as those of Corollary \ref{shacor}. Moreover, over a sufficiently small $U$ these are by construction compatible with the above definition of the pairing (\ref{finalpair}).
\end{dem}

\begin{prop}
For $M_k=[0\to A_k]$ the above pairing reduces to the classical Cassels-Tate pairing.
\end{prop}

\begin{dem}
This was basically proven in \cite{poonen}. There one finds four equivalent definitions of the Cassels-Tate pairing, two of them using the diminished cup-product construction for explicitly defined pairings $A_k(\kbar)\times A_k^{\ast}(\kbar)\to\G[1]$. For instance, one of them (called the Albanese-Picard pairing by \cite{poonen}) is given as follows: replace $A_k(\kbar)$ with the quasi-isomorphic complex $C_1=[Y(A)\to Z(A)]$, where $Z(A)$ is the group of zero-cycles on $A_{\kbar}=A_k\times_k\kbar$ and $Y(A)$ is the kernel of the natural summation map $Z(A)\to A(\kbar)$, and replace $A^{\ast}_k(\kbar)$ by the complex $C_2=[{\kbar (A_{\kbar})}^{\times}/\kbar^{\times}\to {\rm Div}^0(A_{\kbar})]$, where ${\rm Div}^0$ stands for divisors algebraically equivalent to 0. Now define a map $C_1\otimes C_2\to \G(\kbar)[1]$ using the partially defined pairings 
$$
Z(A)\times [{\kbar (A_{\kbar})}^{\times}/\kbar^{\times}]\to {\kbar}^{\times}
$$
and
$$
{\rm Div}^0(A_{\kbar})\times Y(A)\to{\kbar}^{\times}
$$   
where the first is defined on (\cite{poonen}, p. 1116) using the Poincar\'e divisor on $A\times A^{\ast}$ and the second by evaluation. The resulting pairing $A_k(\kbar)\times A_k^{\ast}(\kbar)\to\G[1]$ is one of the classical definitions of the duality pairing for abelian varieties over an algebraically closed field.  
\end{dem}

\begin{rema}\rm Similarly, one verifies that in the case $M_k=[0\to T_k]$ one gets the usual pairing between the Tate-Shafarevich groups of the torus $T_k$ and its character module.
\end{rema}

\begin{rema}\rm
In \cite{flach}, Flach defines a Cassels-Tate pairing for finite dimensional continuous $\ell$-adic representations of the Galois group of a number field $k$ that stabilize some lattice. Comparing the above construction with his shows that in the case when the representation comes from the $\ell$-adic realisation of a 1-motive $M_k$ the two pairings are compatible. In the special cases $M_k=[\Z\to 0]$ and $M_k=[0\to A_k]$ this is already pointed out in Flach's paper (see also \cite{fontaine} for the latter case). So our result can be interpreted as a ``motivic version'' of Flach's pairing for motives of dimension one.
\end{rema}

\bigskip

\noindent {\bf \large Appendix : Completion
of topological abelian groups}

\bigskip

In this appendix we collect some (probably well-known) results that we
needed in the paper.

\begin{propbis} \label{topdual}
Let $0 \ra A \stackrel{i}{\ra} B \stackrel{p}{\ra} C \ra 0$ 
be an exact sequence in the category of topological abelian groups. 

\begin{enumerate}

\item If $i$ is strict with open image, then the map $B^D \ra A^D$ is surjective.

\item Assume that the map $p : B \ra C$ is open. Then the sequences
$$A^{\wedge} \ra B^{\wedge} \ra C^{\wedge} \ra 0$$
$$0 \ra C^D \ra B^D \ra A^D$$
are exact. 

\item Assume that $p$ is open and $i$ is strict with closed image. 
Suppose further that $B$ is Hausdorff, locally compact, 
completely disconnected, and compactly generated. Then the sequences
$$0 \ra A^{\wedge} \ra B^{\wedge} \ra C^{\wedge} \ra 0$$
$$ 0 \ra C^D \ra B^D \ra A^D \ra 0$$
are exact, the map $B \ra B^{\wedge}$ is injective, and $(B^{\wedge})^D=B^D$.

\end{enumerate}

\end{propbis}

\dem 1. If $i$ is strict and $i(A)$ is open in $B$,
then any continuous homomorphism $A \ra \Q/\Z$ induces a continuous
homomorphism $s:\,i(A) \ra \Q/\Z$ which in particular has open kernel. By divisibility of $\Q/\Z$ this extends to a homomorphism
$\bar s : B \ra \Q/\Z$. Since $i(A)$ is open
in $B$, ${\rm Ker} \,\bar s$ is also open in $B$ being a union of cosets of ${\rm Ker}\, s$, whence the continuity  of $\bar s$.

\smallskip

2. To get the first exact sequence, 
the only nontrivial point consists of proving that an element 
$b \in B^{\wedge}$ whose image in $C^{\wedge}$ is trivial comes from 
$A^{\wedge}$. Let $B' \subset B$ be an open subgroup of finite index.
Then $p(B') \subset C$ is of finite index and is open because $p$ 
is open. Since $i$ is continuous, $(A\cap B') \subset A $ is open
and of finite index as well.
Now we use the exact sequence
$$0 \ra A/(A \cap B') \ra B/B' \ra C/p(B') \ra 0$$
where $B'$ runs over the finite index open subgroups of $B$. 

The second exact sequence follows immediately from the fact that the group
$C \simeq B/i(A)$ is equipped with the quotient topology ($p$ 
being an open map). 

\smallskip

3. Let us show the left exactness of the first sequence. 
Since $i$ is strict with closed image, we can assume that $A$ is a closed subgroup of $B$ 
with the induced topology. We have to show that if $A' \subset A$ is 
open and of finite index, then there exists an open subgroup of finite 
index $B' \subset B$ such that $B' \cap A \subset A'$. Replacing 
$A$ and $B$ by $A/A'$ and $B/A'$ (equipped with the quotient topology), 
we reduce to the case when $A$ is a finite subgroup and we must show that 
$B$ contains a finite index open subgroup $B'$ with $B' \cap A=\{ 0 \}$. 
To do so, it is sufficient (using the finiteness of $A$) 
to prove that the intersection of all
finite index open subgroups of $B$ is zero, that is 
$B \hookrightarrow B^{\wedge}$.
By \cite{hewitt}, \Romannumeral2.9.8. (in the special case when $B$ is 
completely disconnected), $B$ is topologically isomorphic 
to a product $\Z^b \times K$ with $K$ compact (hence profinite), thus
$B \hookrightarrow B^{\wedge}=\widehat Z ^b \times K$. 

This also shows that $(B^{\wedge})^D=B^D$. The projection $L$ 
of $A$ on $\Z^b \subset B$ is a discrete lattice and
$A$ is topologically isomorphic to $L \times (A \cap K)$,
hence $(A^{\wedge})^D=A^D$. This proves the exactness of the second sequence.
\enddem

\bigskip
\begin{tabbing}
 \hspace{1cm}\= David Harari\hspace{5cm}\= Tam\a'as Szamuely \\
 \>DMA \> Alfr\a'ed R\a'enyi Institute of Mathematics\\
 \>Ecole Normale Sup\a'erieure \> Hungarian Academy of Sciences\\
  \>45, rue d'Ulm \> PO Box 127 \\ 
 \>F-75230 Paris Cedex 05 \> H-1364 Budapest   \\ 
 \>France \> Hungary \\
 \>{\tt David.Harari@ens.fr} \> {\tt szamuely@renyi.hu}     
 \end{tabbing}

\end{document}